\documentclass [11pt]{article}
\usepackage{graphicx}
\begin{document}

\title{Resolutions of homogeneous bundles on ${\bf P}^2$}
\author{Giorgio Ottaviani \hspace{1cm} Elena Rubei}
\date{}
\maketitle

\def\thefootnote{}
\footnotetext{ \hspace*{-0.36cm}
{\bf Address (of both authors)}: Dipartimento di
Matematica ``U.Dini'', Viale Morgagni 67/A, c.a.p. \hspace*{0.2cm}
 50134 Firenze,
Italia. {\bf E-mail addresses:} ottavian@math.unifi.it,
rubei@math.unifi.it


 }

\newtheorem{theorem}{Theorem}
\newtheorem{lemma}[theorem]{Lemma}
\newtheorem{prop}[theorem]{Proposition}
\newtheorem{rem}[theorem]{Remark}
\newtheorem{cor}[theorem]{Corollary}
\newtheorem{defin}[theorem]{Definition}
\newtheorem{notat}[theorem]{Notation}
\newtheorem{noterec}[theorem]{Notation and recalls}
\newtheorem{exam}[theorem]{Example}






\section{Introduction}

Homogeneous bundles on  ${\bf P}^2=SL(3)/P$ can be described
by representations   of the parabolic subgroup $P$.
In 1966  Ramanan proved  that  if $\rho$ is an irreducible  representation 
of $P$ then the induced bundle $E_{\rho}$ on  ${\bf P}^2 $
is simple and   even stable (see \cite{R}). 
  Since $P$ is not a reductive group, there is a lot of indecomposable
reducible representations of $P$ and
 to classify homogeneous bundles on ${\bf P}^2$ and among them the simple
 ones, the stable ones, etc. 
by means of the study of the representations of the 
parabolic subgroup  $P$ seems difficult.

In this paper our point of view is to consider the minimal free 
 resolution of the bundle. Our aim is to classify 
 homogeneous vector bundles on ${\bf P}^2$ by means of their minimal 
resolutions. 
Precisely we observe that  if $E$ is a
homogeneous vector bundle on ${\bf P}^2={\bf P}(V)$ ($V$ complex
vector space of dimension $3$) there exists 
a minimal free resolution of $E$  $$ 0 \rightarrow \oplus_q
{\cal O}(-q) \otimes_{\bf C} A_{q} \rightarrow \oplus_q {\cal
O}(-q) \otimes_{\bf C} B_{q} \rightarrow E \rightarrow 0 $$ 
with $SL(V)$-invariant maps 
($A_{q}$ and $B_{q}$ are $SL(V)$-representations)
 and we characterize completely the
representations that can occur as $A_{q}$ and $B_{q}$ and the maps
$A \rightarrow B$ that can occur ($A:=\oplus_q {\cal O}(-q)
\otimes_{\bf C} A_{q}$, $ B=\oplus_q {\cal O}(-q) \otimes_{\bf C}
B_{q}$).  To state the theorem we need some notation.

\begin{notat} \label{M(alfa)}
Let $q,r \in {\bf N}$;
for every $\rho \geq p$, let $\varphi_{\rho, p}$ be a fixed  $SL(V)$-invariant
nonzero  map  $S^{p,q,r} V  \otimes {\cal O}_{{\bf P}(V)} (p)
 \rightarrow S^{\rho,q,r} V  \otimes {\cal O}_{{\bf P}(V)} (\rho) $   
(it is unique up to multiples) s.t.
 $\varphi_{\rho , p} = \varphi_{\rho ,p'} \circ \varphi_{p',p}$
$\forall \rho \geq p' \geq p$  (where $S^{p,q,r} $ denotes
 the Schur functor associated to $(p,q,r)$, see  \S 2).

Let ${\cal P}, {\cal R} \subset {\bf N}$, $c\in {\bf Z}$;
 for any $SL(V)$-invariant map $$\gamma:
 \oplus_{p \in {\cal P}}  A^{p} \otimes S^{p,q,r}V (c+p) 
  \rightarrow 
 \oplus_{\rho  \in {\cal R}}  B^{\rho} \otimes S^{\rho,q,r}V (c+\rho) $$
($ A^{p}$ and $ B^{\rho}$ vector spaces)
 we define    $M(\gamma)$ to be the map
$$M(\gamma) : \oplus_{p \in {\cal P}}  A^{p} \rightarrow 
 \oplus_{\rho \in {\cal R}}  B^{\rho}$$ s.t.
$\gamma_{\rho, p} = M(\gamma)_{\rho, p} \otimes \varphi_{\rho,p}$ $\forall
\rho, p$  with $\rho \geq p$ and $ M(\gamma)_{\rho, p}=0 $  $\forall
\rho, p $ with $\rho < p$,  where $\gamma_{\rho, p} : 
A^{p} \otimes  S^{p,q,r} V (c+p)  \rightarrow  
B^{\rho} \otimes   S^{\rho,q,r}V (c+\rho)  $ and
$ M(\gamma)_{\rho, p}:  A^{p} \rightarrow   B^{\rho}$
are the maps induced respectively by $\gamma$ and $M(\gamma)$
(by restricting and projecting).
\end{notat}

\begin{theorem} \label{fibomintro}
i) On ${\bf P}^2= {\bf P}(V)$
 let $A= \oplus_{p,q,i}  A_i^{p,q} \otimes S^{p,q}V (i)$, $
 B= \oplus_{p,q,i} B_i^{p,q} \otimes S^{p,q}V (i) $ with 
$p,q,i$ varying in  a finite  subset of ${\bf N}$,   $A_i^{p,q}$ and 
 $B_i^{p,q}$ finite dimensional vector spaces. 
Then $A$ and $B$  are the first two terms of a minimal free 
resolution of a homogeneous bundle on ${\bf P}^2$  if and only if   
$\forall c \in {\bf Z}, q, \tilde{p} \in {\bf N}$  $$dim(\oplus_{p \geq
\tilde{p}} A_{c+ p}^{p,q}) \leq dim (\oplus_{\rho >
\tilde{p}} B_{c+\rho}^{\rho,q})$$

ii)  Let $$A= \oplus_{p,q,r} A^{p,q,r} \otimes S^{p,q,r}V (p+q+r) 
\hspace{1,3cm} B= \oplus_{p,q,r} B^{p,q,r} \otimes S^{p,q,r}V (p+q+r)
$$ $p,q,r$ varying in a finite subset of ${\bf N}$,
$A^{p,q,r}$ and $B^{p,q,r}$ finite dimensional vector spaces;
 let $\alpha $ be  an $SL(V)$-invariant map
$A \rightarrow B$.
Then there exists a homogeneous bundle $E$ on ${\bf P}^2$  s.t.
$$ 0 \rightarrow A \stackrel{\alpha}{\rightarrow} B \rightarrow E 
\rightarrow 0$$ is a minimal free resolution of $E$ if and only if 
$ M(\alpha_{p,q,r}) :
  A^{p,q,r}  \rightarrow   B^{p,q,r} 
 $ is zero $\forall p,q,r$ and 
$M(\alpha_{q,r}) :\oplus_{p}  A^{p,q,r} \rightarrow \oplus_{\rho}
  B^{\rho,q,r}$ is injective  $\forall q,r$, where  
$ \alpha_{p,q,r} :
  A^{p,q,r} \otimes S^{p,q,r} V (p) \rightarrow   B^{p,q,r} \otimes S^{p,q,r}
 V (p)$ and  $\alpha_{q,r} :
\oplus_{p}  A^{p,q,r} \otimes S^{p,q,r} V (p)  \rightarrow  \oplus_{\rho}
 B^{\rho,q,r} \otimes
 S^{\rho,q,r} V (p) $  are 
the maps induced by $\alpha$. 
\end{theorem}

The above theorem   allows us to parametrize the set of
 homogeneous bundles on ${\bf P}^2$ by a set of sequences
 of injective matrices 
with a certain shape up to the action of invertible matrices with a certain 
shape. 
An interesting  problem is to use this 
parametrization to study the related  moduli spaces.

Then    we begin to  study which minimal free resolutions
give simple or stable homogeneous bundles.  
We will consider the case $A$ is irreducible;  we call a homogeneous bundle
{\it elementary}  if it has  minimal free resolution $0 \rightarrow
 A \rightarrow B \rightarrow E \rightarrow 0$
 with $A$ irreducible; besides 
we say that a bundle $E$ on ${\bf P}^2$ is {\it regular}  if the minimal free 
resolution is $0 \rightarrow A \rightarrow B \rightarrow E \rightarrow 0$ with 
 all the components of $A$ with the same twist 
 and all the components of $B$ with the same twist.

First we study simplicity and stability of regular elementary
 homogeneous bundles.
Fundamental tool are quivers and representations of  quivers associated to
 homogeneous bundles introduced by Bondal and Kapranov in \cite{B-K}.
The quivers allow us to handle well and ``to make explicit'' the homogeneous 
subbundles of a homogeneous bundle $E$
and Rohmfeld's criterion (see \cite{Rohm}) 
 in this context is equivalent to saying that 
$E$ is semistable if and only if the slope   of every subbundle associated 
to a subrepresentation of the quiver representation of $E$ is less or equal 
than the slope of $E$.

The simplest regular elementary  homogeneous bundles are the 
bundles $E$ defined by the exact sequence
$$ 0 \rightarrow S^{p,q} V \otimes
{\cal O}(-s)  \stackrel{\varphi}{\rightarrow}  S^{p+s,q} V \otimes
{\cal O}  \rightarrow E \rightarrow 0 $$
for some $p, q ,s \in {\bf N}$, $p\geq q$, 
 $\varphi$  an $SL(V)$-invariant nonzero map; 
we prove that such bundles are stable, see Theorem \ref{stable}
(observe that Ramanan's theorem  does not apply here). 
In the particular case   $p=q=0$
the stability  of $E$ was already proved    in \cite{Ba}.

Besides we prove

\begin{theorem} \label{pinco}
A  regular elementary  homogeneous bundle $E$ on ${\bf P}^2$
 is simple if and only if its minimal free 
resolution is of the following kind: $$  0 \rightarrow S^{p,q} V \otimes
{\cal O}(-s)   \stackrel{\varphi}{ \rightarrow} W \otimes
{\cal O}  \rightarrow E \rightarrow 0  $$
where  $p, q ,s \in {\bf N}$, $p\geq q$, 
 $W $ is a nonzero $SL(V)$-submodule of $S^{p,q} V \otimes S^s V$,
all the components of $\varphi $ 
are nonzero $SL(V)$-invariant maps 
and we are in one of the following cases:

i) $p=0$ 

ii) $p> 0$ and
 $W  \neq S^{p,q} V \otimes S^s V$.
\end{theorem}

By using the above theorem and  Theorem \ref{1regular}, which characterizes 
stability of the minimal  free resolution of regular elementary homogeneous 
bundles when (with the above notation) $s=1$, we find  infinite examples 
of unstable simple  homogeneous bundles.

Finally  we state a criterion, generalizing Theorem \ref{pinco}, 
 to say when an elementary (not necessarily 
regular) homogeneous bundle is simple by means of its minimal free resolution,
see Theorem \ref{fibomsem}.

The sketch of the paper is the following: in \S 2 we recall some basic facts
on representation theory; in \S 3  we characterize the resolutions of 
homogeneous bundles on ${\bf P}^2$: 
 in this section we prove Theorem \ref{fibom}, which contains Theorem 
\ref{fibomintro}; in \S 4
 we recall the theory of quivers; in \S 5  we prove some lemmas by 
using  quivers and we fix some notation;   in \S 6 we study stability 
and simplicity of   elementary  homogeneous bundles.

\section{Notation and preliminaries}

We recall some facts  on representation theory (see for instance \cite{F-H}).

 Let $d$ be a natural number and let $ \lambda =
(\lambda_1,...,\lambda_k)$ be a partition of $d$ with $\lambda_1 \geq ...
\geq \lambda_k$.

We can associate to $ \lambda $ a {\bf Young diagram} with
$\lambda_i$ boxes in the $i$-th row, the rows lined up on the
left.
The conjugate partition $\lambda'$ is the partition of $d$
whose Young diagram is obtained from the Young diagram of $\lambda$
interchanging rows and columns.

A {\bf tableau} with entries in $\{1,..., n\}$
 on the Young diagram of a partition  $\lambda
=(\lambda_1,...,\lambda_k)$ of $d$ is a numbering of the boxes by
the integers $ 1,...., n$, allowing repetitions 
(we say also that it is a tableau on $\lambda$).

\begin{center}
\includegraphics[scale=0.38]{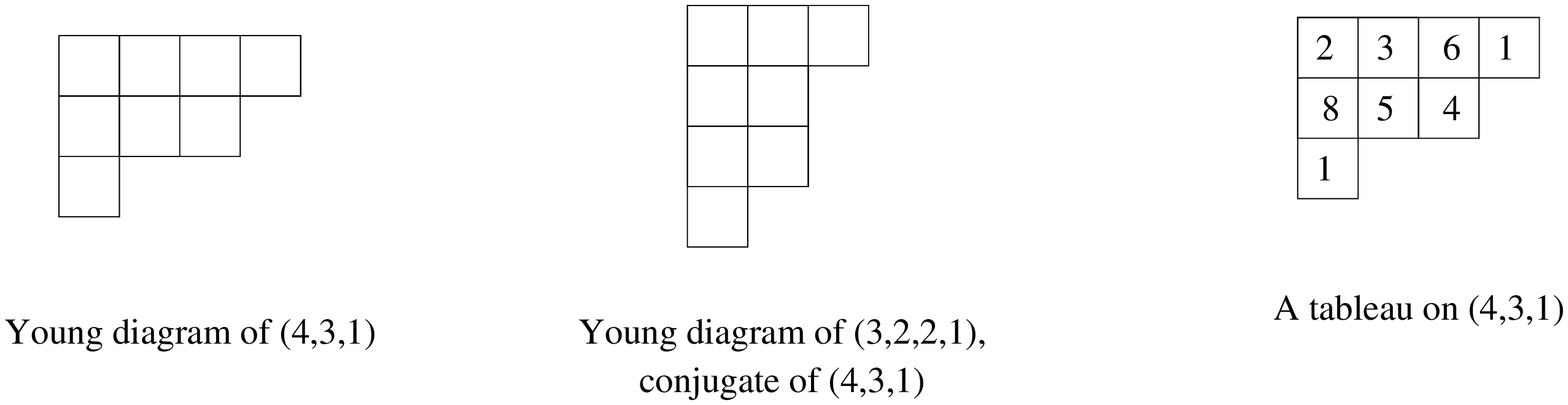}
\end{center}

\begin{defin}
Let $V$ be a complex vector space of dimension $n$. 
Let $d \in {\bf N}$  and let $ \lambda =
(\lambda_1,...,\lambda_k)$ be a  partition of $d$, with $\lambda_1 \geq ...
\geq  \lambda_k$.
Number the boxes of the Young diagram of $\lambda$ with the numbers $1,...,d$
from left to right beginning from the top row. 
Let $\Sigma_d$ be the group of permutations on $d$ elements; 
let $R$ be the subgroup of $\Sigma_d$ given by the permutations preserving
the rows and 
let $C$ be the subgroup of $\Sigma_d$ given by the permutations preserving
the columns.

We define $$S^{\lambda} V := Im (\sum_{a \in C , s \in R} sign(a) s \circ a:
\otimes^d V \rightarrow \otimes^d V )$$
The $S^{\lambda} V$ are called {\bf Schur representations}.
\end{defin}


The $S^{\lambda} V $ are irreducible  $SL(V)$-representations and it is
well-known that all the irreducible $SL(V)$-representations are of this form.

\begin{notat} \label{simant}
$ \bullet  $ Let $V$ be a complex vector space 
and let $\{v_j\}_{1,..., n}$
be a basis of $V$. Let $d \in {\bf N} $ and 
$\lambda =(\lambda_1, ...,\lambda_k)$ be a
partition of $d$ with $\lambda_1 \geq ... \geq \lambda_k$.
Let $ \mu $ be the conjugate partition. Let ${\cal G}_{\lambda}$ be 
 the free abelian group generated  by
the tableaux on $\lambda$ with entries in  $\{1,...,n\} $
and let ${\cal T}_{\lambda} = {\cal G}_{\lambda} 
\otimes_{{\bf Z}}
 {\bf C} $.
 
$\bullet$ 
Let $t$ be the map associating to a tableau $T$ on $ \lambda$ with entries
 in $\{1,..., n\}$
the following element of $V^{\otimes d}$: $$ v_{T_1^1} \otimes....
   \otimes  v_{T_{\lambda_1}^1} \otimes ....\otimes v_{T_1^k}
    \otimes....
   \otimes  v_{T_{\lambda_k}^k}$$ where $(T_1^j ... T_{\lambda_j}^j)$
   is the $j$-th   row of $T$.

$\bullet$   We define $ant:  {\cal T}_{\lambda} \rightarrow 
{\cal T}_{\lambda} $ to be the linear map 
   s.t. for every  tableau $T$ on $\lambda$  $$ant(T) =\sum_{(\sigma_1 , ...,
\sigma_{\lambda_1}) \in \Sigma_{\mu_1} \times ....\times
\Sigma_{\mu_{\lambda_1}} } sign (\sigma_1)  ...
sign(\sigma_{\lambda_1}) T^{\sigma_1, ...,
\sigma_{\lambda_1}}$$ where 
$T^{\sigma_1, ..., \sigma_{\lambda_1}}$ is the tableau obtained
from $T$ permuting the elements of the $j$-th column with
$\sigma_j$ $\forall j$. For instance

\begin{center}
\includegraphics[scale=0.42]{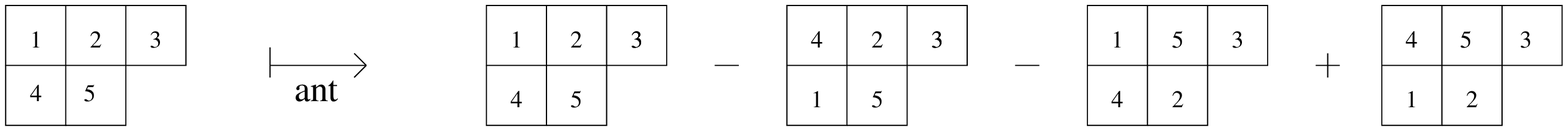}

\end{center}

$\bullet$  Analogously we define $sim : 
 {\cal T}_{\lambda} \rightarrow {\cal T}_{\lambda} $ to be the linear map
 s.t. for every 
tableau $T$ on $\lambda$  $$sim(T)= \sum_{(\sigma_1 , ..., \sigma_{\mu_1}) \in
\Sigma_{\lambda_1} \times ....\times \Sigma_{\lambda_{\mu_1}} }
T_{\sigma_1, ..., \sigma_{\mu_1}}$$ where $T_{\sigma_1, ...,
\sigma_{\mu_1}}$ is the tableau obtained from $T$ permuting the
elements of the $j$-th row with $\sigma_j$ $\forall j$. For instance

\begin{center}
\includegraphics[scale=0.4]{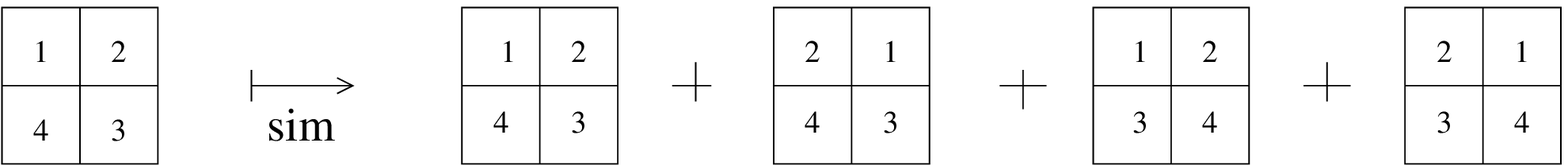}
\end{center}

$\bullet$ We call $ord $ the linear map  associating to a tableau
$T$ the tableau obtained from $T$ ordering the entries of every
row in  nondecreasing way. 

Observe that $ ord \circ sim = \lambda_1 ! ... 
\lambda_k ! \,  ord $ and
$ sim \circ ord = sim$.

$\bullet $ Let ${\cal S} = t \circ sim \circ ant$.

\end{notat}


Obviously 
the space $S^{\lambda}V$ can be described as the image of
 ${\cal S} : {\cal T}_{\lambda}
 \rightarrow V^{\otimes d} $.

\bigskip

We recall that {\bf Pieri's formula} says that, 
 if $\lambda=(\lambda_1, \lambda_2,...) $ is a partition of a natural
number  $d$ with $\lambda_1 \geq \lambda_2 \geq ..$
and $t$ is a natural number, then $$S^{\lambda} V \otimes 
S^t V = \oplus_{\nu } S^{\nu} V$$  as $SL(V)$-representations,
where the sum is performed on all the 
partitions $\nu= (\nu_1,..) $ with $\nu_1 \geq \nu_2\geq ...$ 
 of $d+t$ whose Young diagrams are
obtained from the Young diagram of $\lambda$ adding  $t$ boxes
not two in the same column.

Finally we observe  that 
if $V$ is a complex vector space of
dimension $n$ then $S^{(\lambda_1,...,\lambda_{n-1})} V$
  is isomorphic to $S^{(\lambda_1+r,...,\lambda_{n-1} +r,r) }V$  for all $r$
as $SL(V)$-representation. Besides
$(S^{(\lambda_1,...,\lambda_{n})}V)^{\vee}$ is isomorphic as
$SL(V)$-representation to $S^{ (\lambda_1 -\lambda_n,..., 
\lambda_1-\lambda_{2}) }V$.
Moreover $(S^{\lambda}V)^{\vee} \simeq S^{\lambda} V^{\vee}$.

\begin{notat}
$\bullet $ 
In all the paper  $V$ will be a complex vector space of dimension $3$ if 
not otherwise specified.

$\bullet$ If $E$ is a vector bundle on ${\bf P}(V) $ then $\mu(E)$ will 
denote the slope of $E$, i.e. the first Chern class divided by the rank.
\end{notat}

\section{Resolutions of homogeneous vector bundles}

The aim of this section is to characterize the minimal free resolutions of 
the homogeneous bundles on ${\bf P}^2$. 

\begin{lemma} Let $E$ be a homogeneous vector bundle on ${\bf
P}^2={\bf P}(V)$.
By Horrocks' theorem \cite{Hor} the bundle $E$ has a minimal free 
resolution
 $$ 0 \rightarrow \oplus_q {\cal O}(-q) \otimes_{\bf C} A_{q}
\rightarrow \oplus_q {\cal O}(-q) \otimes_{\bf C} B_{q}
\stackrel{\psi}{\rightarrow} E \rightarrow 0 $$ 
 ($A_q$, $B_q$ ${\bf C}$-vector spaces). 
Since $E$ is homogeneous we can suppose the maps are $SL(V)$-invariant
 maps ($A_{q}$ and $B_{q}$ are $SL(V)$-representations).
\end{lemma}

{\it Proof.}   (See also \cite{Ka}).
Let  $M:= \oplus_k H^0 (E(k))$ and $S = \oplus_k Sym^k V^{\vee} $. 
Let $m_i \in M$  $i=1,...,k$ s.t. no one of them can be written 
as linear 
 combination with coefficient in $S$ of elements of the $SL(V)$-orbits of 
the others and s.t. $\cup_{g \in SL(V), i=1,..,k}
  g (m_i)  $ generates $M$  on $S$. Let $q_i =deg (m_i)$.
Let $B_i :=  \langle \cup_{g \in SL(V)}  g (m_i)  
\rangle_{{\bf C}}$ (finite dimensional $SL(V)$-representations).
Let $P= \oplus_i B_i \otimes_{{\bf C}} S(-q_i) $ and $\varphi : P 
\rightarrow M$ be the $SL(V) $-invariant map given by multiplication.
Let 
$\psi : B \rightarrow E $ be the sheafification of 
$\varphi: P \rightarrow M$.
 Thus $B =\oplus_i B_i \otimes {\cal O}(-q_i) $. Let 
$A= Ker (\psi)$; it is a homogeneous vector bundle; we have $H^1 (A(t) ) =0 $
 $\forall t \in {\bf Z} $, because $ H^0 (B(t)) \rightarrow H^0(E(t))  $ is 
surjective      $\forall t $ and $H^1 (B(t)) = 0$  $\forall t$, hence by 
Horrocks'  criterion $A$ splits.
\hfill  \framebox(7,7)

\begin{rem} {\rm 
a) If $U,W,V$ are three  vector spaces then  on $P(V)$ we have 
$Hom (U \otimes {\cal O} (-s)  ,W) = Hom (U \otimes S^s V   ,W) $ (the 
isomorphism can be given by $H^0 (\cdot^{\vee})^{\vee} $).

b) Let $V$ be a vector space. For any $\lambda$, $\mu $ partitions, 
$s \in {\bf N}$, up to multiples there is a unique
 $SL(V)$-invariant map  
$$ S^{\lambda} V \otimes {\cal O}(-s) \rightarrow S^{\mu} V \otimes {\cal O}$$
by Pieri's formula,  Schur's lemma and part a of the remark. }
\end{rem}

Theorem \ref{fibom}, which implies Theorem \ref{fibomintro}, 
 is the aim of this section. It allows us
to classify all homogeneous vector bundles  on ${\bf P}^2$; in
fact it  characterizes their minimal free  resolutions.

Precisely part {\it i} allows us to say which $A$ and $B$ can occur in
a minimal free resolution  $0 \rightarrow A \rightarrow B \rightarrow E
\rightarrow 0$ of a homogeneous vector bundle $E$ on ${\bf P}^2$; first we
 investigate   when, given $A$ and $B$ 
direct sums of bundles of the kind $S^{p,q}V (-i) $,  
 there exists an injective $SL(V)$-invariant map $A \rightarrow B$; roughly 
speaking this is true 
 if and only if  for every $SL(V)$-irreducible subbundle  $S=  S^{p,q} V (i) $ 
 of $A$  there  exists  a subbundle $M(S)$
 of $B$ of the kind $ S^{p+s,q} V (i+s) $
for some $s \in {\bf N}$ and we can choose $M(S)$ in such way that
 the map $S \mapsto M(S)$ is  injective. A crucial point of  the proof 
is the fact that an $SL(V)$-invariant map $S^{p,q} V 
\rightarrow S^{p+s_1,q +s_2,s_3 }V (s_1 +s_2 +s_3 ) $ 
is injective if and only if $s_2=s_3=0$
and  the intersection of the kernels of the ones of 
such maps with $s_2 +s_3 >0 $ is 
nonzero.

Part {\it ii} allows us to say  which maps  $A
\rightarrow B$ can occur in a minimal free resolution 
$ 0 \rightarrow A \rightarrow B \rightarrow E \rightarrow 0$ of a homogeneous
 bundle $E$; first  we study  when   an
$SL(V)$-invariant map $\alpha: A \rightarrow B$  is injective ($A$ and $B$
direct sums of bundles of the kind $S^{p,q,r}V (-i) $).
We remark that,  if $\alpha: A \rightarrow B$ is 
 an $SL(V)$-invariant map,
 we can suppose that the sum $p+q+ r-i $ is constant 
 for  $ S^{p,q,r}V (i) $ varying among all $SL(V)$-irreducible subbundles 
of $A$ or $B$ (by using the isomorphism $S^{p,q} V
\simeq S^{p + u,q+u,u}V$  $\forall u \in {\bf N}$).

In the sequel we will use Notation \ref{M(alfa)}.

\begin{theorem} \label{fibom}
i) On ${\bf P}^2= {\bf P}(V)$
 let $A= \oplus_{p,q,i}  A_i^{p,q} \otimes S^{p,q}V (i)$, 
$ B= \oplus_{p,q,i} B_i^{p,q} \otimes S^{p,q}V (i) $ with 
$p,q,i$ varying in  a finite  subset of ${\bf N}$,   $A_i^{p,q}$ and 
 $B_i^{p,q}$ finite dimensional vector spaces. 
There
exists an injective $SL(V)$-invariant map $A \rightarrow B$ if and
only if $\forall c \in {\bf Z}, q, \tilde{p} \in {\bf N}$
\begin{equation} \label{grr}  dim(\oplus_{p \geq
\tilde{p}} A_{c+ p}^{p,q}) \leq dim (\oplus_{\rho \geq
\tilde{p}} B_{c+\rho}^{\rho,q})\end{equation}
i') Besides $A$ and $B$  are the first two terms of a minimal free 
resolution of a homogeneous bundle on ${\bf P}^2$  if and only if 
$\forall c \in {\bf Z}, q, \tilde{p} \in {\bf N}$  
 \begin{equation} \label{zzz} dim(\oplus_{p \geq
\tilde{p}} A_{c+ p}^{p,q}) \leq dim (\oplus_{\rho >
\tilde{p}} B_{c+\rho}^{\rho,q}) \end{equation}
ii)  Let $$A= \oplus_{p,q,r} A^{p,q,r} \otimes S^{p,q,r}V (p+q+r) 
\hspace{1,3cm} B= \oplus_{p,q,r} B^{p,q,r} \otimes S^{p,q,r}V (p+q+r)
$$ $p,q,r$ varying in a finite subset of ${\bf N}$,
$A^{p,q,r}$ and $B^{p,q,r}$ finite dimensional vector spaces;
 let $\alpha $ be  an $SL(V)$-invariant map
$A \rightarrow B$; $\alpha $ is injective if and only if $\forall
q,r $ the induced map
 $$\alpha_{q,r} :\oplus_{p }
A^{p,q,r} \otimes S^{p,q,r}V (p+q+r)  \rightarrow \oplus_{\rho
} B^{\rho,q,r} \otimes S^{\rho,q,r}V (\rho+q+r)$$ 
is injective and (by Lemma \ref{prelim})  this is true if and only if 
$M(\alpha_{q,r}) :\oplus_{p}  A^{p,q,r} \rightarrow \oplus_{\rho}
  B^{\rho,q,r}$ is injective.

Besides obviously there exists a homogeneous bundle $E$ on ${\bf P}^2$  s.t.
$$ 0 \rightarrow A \stackrel{\alpha}{\rightarrow} B \rightarrow E 
\rightarrow 0$$ is the minimal free resolution of $E$ if and only if 
$M(\alpha_{q,r}) :\oplus_{p}  A^{p,q,r} \rightarrow \oplus_{\rho}
  B^{\rho,q,r}$ is injective  $\forall q,r$ and $M(\alpha_{p,q,r}) :
  A^{p,q,r} \rightarrow   B^{p,q,r}$ is zero $\forall p,q,r$. 
\end{theorem}

To prove Theorem  \ref{fibom} we need some lemmas.

\begin{lemma} \label{iniet} Let $V$ be a complex vector space of dimension 
$n$. Let $\lambda_1,...., \lambda_{n-1}, s \in {\bf N} $ with 
$ \lambda_1 \geq ...  \geq \lambda_{n-1} $ and
$$\pi: S^{\lambda_1,..., \lambda_{n-1}} V
  \otimes S^s V \rightarrow S^{\lambda_1+s,\lambda_2, ..., \lambda_{n-1}} V$$
 be  an
$SL(V)$-invariant nonzero map. Then (up to multiples)
$\pi$ can be described in the following way:  let $T$ be a tableau
on $(\lambda_1,..., \lambda_{n-1})$
 and let $R$ be a tableau on $(s)$; then $$\pi({\cal S}(T) \otimes
{\cal S}(R))={\cal S}(TR)$$ 
 where $TR$ is the tableau obtained from $T$ adding
$R$ at the end of  its first row and ${\cal S} $ is defined in Notation 
\ref{simant}.
In particular  on ${\bf P}^{n-1} ={\bf P}(V) $  any $SL(V)$-invariant
nonzero map $$ S^{\lambda_1,..., \lambda_{n-1}} 
V (-s) \rightarrow S^{\lambda_1 +s , \lambda_2 ,..., \lambda_{n-1}} V$$ 
 is injective.
\end{lemma}

{\it Proof.} Let $\varphi: S^{\lambda_1,..., \lambda_{n-1}} V
  \otimes S^s V \rightarrow S^{\lambda_1+s,\lambda_2, ..., \lambda_{n-1}} V$ 
be the linear map s.t.
 $\varphi ({\cal S}(T) \otimes
{\cal S}(R)) = {\cal S}(TR) $  $\forall T$  tableau on $(\lambda_1,..., 
\lambda_{n-1})$ and $\forall R$  tableau on $(s)$; 
it is sufficient to prove that $\varphi$ is well defined, $SL(V)$-invariant
and nonzero.

To show that it is well defined it is sufficient to show that, if
$T, T' \in {\cal T}_{(\lambda_1,..., \lambda_{n-1})}$
s.t. ${\cal S}(T) = {\cal S}(T')$ and $R,R' \in {\cal T}_{(s)}$ 
s.t. ${\cal S}(R)={\cal S}(R')$,
   then $ {\cal S}(TR) ={\cal S}(T'R')$ (with the obvious definition of $TR$).
Observe that {\small $$  {\cal S}(TR) ={\cal S}(T'R') \Leftrightarrow 
ord (ant (TR))=ord ( ant(T'R')) \Leftrightarrow ord ( ant (T) R) =ord ( ant
(T')R')$$} \hspace*{-0.356cm} and the last  equality  follows from the fact
 that $ord(ant(T))
=ord(ant(T'))$ because $sim(ant(T)) =sim(ant(T'))$ and $ord(R) =
ord(R')$ because $sim(R)=sim(R')$.

Besides obviously the map
 ${\cal S}(T) \otimes {\cal S}(R) \mapsto {\cal S}(TR) $ is
$SL(V)$-invariant and nonzero.
Thus, up to multiples, it is the map $\pi$.

This implies the injectivity of any $SL(V)$-invariant nonzero
 map  $S^{\lambda_1,...,\lambda_{n-1}} 
V (-s) \rightarrow S^{\lambda_1+s, \lambda_2 , ...,\lambda_{n-1} } V$;
 in fact the induced map on 
the fibre on $[0:...:0:1]$ is $$ {\cal S}(T) \mapsto {\cal S}(T [n...n])$$ 
$\forall T \in {\cal T}_{(\lambda_1,..., \lambda_{n-1})}$  and 
if $ {\cal S}(T [n...n]) =0$, then  $ ord \circ sim \circ ant  (T [n...n])=0$;
  thus  $ ord \circ ant  (T [n...n])=0 $, but   
$ord \circ ant  (T [n...n])= 
   (ord  \circ ant  (T)) [n...n]$, hence $ord \circ ant (T) =0$, i.e. 
${\cal S}(T)= 0$. 
 \hfill  \framebox(7,7)

\medskip

The injectivity statement of Lemma \ref{iniet} (probably well known)
will be obvious by the theory of quivers, precisely 
it will follow from Lemma \ref{palle}, but we wanted to show the above
 proof because it is more elementary and intuitive.

\begin{lemma} \label{xazero} {\rm Let  
${\bf P}^{n-1} ={\bf P}(V) $.
 For every $\lambda_1,...,\lambda_{n}
 \in {\bf N}$, $t \in {\bf Z}$  with $\lambda_1 \geq .... \geq \lambda_n$, 
   $ \exists \; y= y_{\lambda_1,..., \lambda_n}^t \in S^{\lambda_1,..., 
\lambda_n}V (t)$,
 $y \neq 0 $ s.t.
 $\varphi (y )= 0$
for every   $SL(V)$-invariant map
$$\varphi :S^{\lambda_1,...., \lambda_n}V(t) 
\rightarrow S^{\lambda_1 +s_1,....,\lambda_n+s_n} V(t+s_1+ ....+s_n)$$
 $\forall \;s_1,...,s_n \in {\bf N}$  s.t. $s_2+...+s_n> 0$.}
\end{lemma}

{\it Proof.} It is sufficient to prove the statement when $t=0$.
It is sufficient to take as $y$ 
 a nonzero element of the image
of an $SL(V)$-invariant nonzero map $$\psi:  S^{\lambda_2,..., \lambda_n}
V (-\lambda_1)\rightarrow
S^{\lambda_1,..., \lambda_n} V  $$ 
(such  a map exists  because, by Pieri's formula,
 $S^{\lambda_1,..., \lambda_n} V $ 
is a summand  of $S^{\lambda_2,..., \lambda_n} V \otimes S^{\lambda_1} V$, 
thus we can take the map  induced by the projection;
besides an $SL(V)$-invariant map $$S^{\lambda_2,..., \lambda_n} V 
(-\lambda_1 -s_1 ...-s_n) \rightarrow  
S^{\lambda_1+s_1, ...., \lambda_n+s_n} V$$  is  zero $\forall s_1,...,s_n$
 with $ s_2 +...+s_n > 0$, 
because by  Pieri's formula the induced map (i.e. $H^0 ( \cdot^{\vee})^{\vee}$)
$ S^{\lambda_2,...,\lambda_n } 
V \otimes S^{\lambda_1+s_1 +... +s_n} V  \rightarrow S^{\lambda_1+s_1, ....,
 \lambda_n +s_n} V$
 is zero, thus $\varphi 
\circ \psi =0$ $\forall \varphi$ as in the statement and then
$\varphi (y)=0 $). 
\hfill  \framebox(7,7)

\begin{rem} \label{piccolo} {\rm 
Let $W$ be a finite dimensional ${\bf C}$-vector
 space 
and $v_1,..., v_k \in W$ not all zero. Let $B$ be a matrix with $k$ columns
s.t. $\sum_{j=1,..., k} 
B_{i,j} v_j =0$ $\forall i$ (i.e. the coefficients of every  row
of $B$ are the coefficients of a linear relation among the $v_j$). Then $B: 
{\bf C}^k \rightarrow {\bf C}^s$ (where $s$ is the number of the rows of $B$) 
is not injective.}
\end{rem}

\begin{lemma} \label{prelim}
Fix  $ q,r \in {\bf N}$ with $q \geq r$ and let ${\cal P}, {\cal R} $ be
 finite subsets of $ \{ p \in {\bf N} | \; p \geq q\}$. 
A map $$\alpha :\oplus_{p \in {\cal P}}
 A^{p} \otimes S^{p,q,r}V (p)   \rightarrow
\oplus_{\rho \in {\cal R}} B^{\rho} \otimes S^{\rho,q,r}V (\rho) 
$$  ($A^p$, $B^p$ nonzero finite-dimensional vector spaces)
is injective if and only if $M(\alpha)$ is injective.
\end{lemma}

{\it Proof.}  Let $\underline{p} = min \; {\cal P}$ and
 $\overline{p} = max \; {\cal P}$.
Let $$\psi=
  \oplus_{p \in {\cal P} }  I_{A^{p}} \otimes
\varphi_{p, \underline{p}} : (\oplus_{p \in {\cal P}, \; p \geq \underline{p}}
A^{p}) \otimes S^{\underline{p},q,r}V
 (\underline{p}) \rightarrow \oplus_{
p \in {\cal P}} A^{p} \otimes S^{p,q,r}V (p)$$
The map $\psi$ 
  is injective by Lemma \ref{iniet}.
 We have \begin{equation} 
\label{alfapsi} \alpha
 \circ \psi =
 (\oplus_{\rho \in {\cal R} , \rho \geq \underline{p}}
I_{B^{\rho}} \otimes \varphi_{\rho , \underline{p}} )\circ 
(M(\alpha) \otimes I_{S^{\underline{p},q,r}V 
(\underline{p})}) \end{equation}

Suppose $\alpha$ is injective.
Since   $\psi $ is injective, $\alpha \circ \psi $ is injective.
Thus, by (\ref{alfapsi}), $M(\alpha)$ is injective.

Suppose now $M(\alpha) $ is injective.
Let $x \in {\bf P}^2$. Let $\alpha^x$ be 
 the map induced on the fibres on $x$ by $\alpha$
and for any bundle $E$, $E^{x}$ will 
denote the fibre on $x$.
 
Let ${\cal P}= \{p_1,..., p_n\}$.  
Let $v  = (v^1_1, ...,v^1_{a_1},......................,v^n_1, ...,v^n_{a_n}
   ) \in Ker ( \alpha^{x}) $  
where $a_i =dim \; A^{p_i}$ and 
$v^{i}_{j} \in S^{p_i,q,r}V (p_i)^x \ $.
We want to show $v=0$.

We can see every $v^{i}_j $ in $S^{\overline{p},q,r} (\overline{p})^x$ (by 
$\varphi_{\overline{p},p_i}$, which is an injection); (after fixing bases 
and seeing $M(\alpha)$ as a matrix)   the coefficients
of the rows of $ M(\alpha)$  are the coefficients of linear relations among the
$v^{i}_j$ seen in $S^{\overline{p},q,r} (\overline{p})^x$, since 
$\alpha^x(v)=0$. Thus, since $M(\alpha) $ is injective,
 by Remark \ref{piccolo}, the $v^{i}_j$ must be all zero, i.e. $v=0$. 
\hfill  \framebox(7,7)

\bigskip

{\it Proof of Theorem \ref{fibom}.} {\it i)} 
Let $ \alpha:A \rightarrow B$ be an $SL(V)$-invariant injective
map.

For any $ c,q, \tilde{p}$, 
let  $$\alpha_{c,q, \geq  \tilde{p}}: \oplus_{p \geq \tilde{p}}
A_{c+p}^{p,q} \otimes S^{p,q}V (c +p)   \rightarrow \oplus_{\rho \geq
\tilde{p}}  B_{c +\rho}^{\rho,q} \otimes S^{\rho,q}V (c +\rho) $$ the
map induced by $\alpha$,  and 
let $$\psi_{\tilde{p}}= 
 \oplus_{p \geq \tilde{p}}  I_{A_{c+p}^{p,q}} \otimes
\varphi_{p, \tilde{p}} (c) : (\oplus_{p \geq \tilde{p}}
A_{c+p}^{p,q}) \otimes  S^{\tilde{p},q}V (c +\tilde{p}) \rightarrow
\oplus_{p \geq \tilde{p}}
 A_{c+ p}^{p,q} \otimes S^{p,q}V (c +p)$$
The map $\psi_{\tilde{p}} $ is injective by Lemma \ref{iniet}.
We can write \begin{equation} \label{ast} \alpha_{c,q, \geq 
\tilde{p}} \circ \psi_{\tilde{p}} =  (\oplus_{\rho \geq \tilde{p}}
I_{B_{c +\rho}^{\rho,q}} \otimes \varphi_{\rho , \tilde{p}} )\circ
(M(\alpha_{c,q, \geq  \tilde{p}}) \otimes 
I_{S^{\tilde{p}, q} V (c + \tilde{p})})  \end{equation}
Let $v \in Ker M(\alpha_{c,q, \geq  \tilde{p}})$. If $v \neq 0$ then 
$\psi_{\tilde{p}} (v \otimes 
y^{c + \tilde{p}}_{\tilde{p},q})$ 
(see Lemma \ref{xazero} for the definition of 
$y^{c + \tilde{p}}_{\tilde{p},q}$) is nonzero (since 
$\psi_{\tilde{p}} $ is injective) and it is  in
 $Ker (\alpha)$ (in fact it is in  $Ker(\alpha_{c,q, \geq  
\tilde{p}})$
  by (\ref{ast})  and thus in $Ker (\alpha) $ by the definition of 
$y^{c + \tilde{p}}_{\tilde{p},q}$). Since $\alpha $ is injective 
we get a contradiction, thus $v =0$.
Thus  $ M(\alpha_{c,q, \geq  \tilde{p}})$ is injective.
Then (\ref{grr}) holds. 

Suppose now (\ref{grr}) holds.
Since $$A=\oplus_{c,q} (\oplus_{p} A_{c+p}^{p,q} \otimes S^{p,q}V (c+p))  
\;\;\;\;\;
B=\oplus_{c,q} (\oplus_{\rho}  B_{c+ \rho}^{\rho,q} \otimes S^{\rho,q}V
(c+\rho) ),$$ to find an injective map
$\alpha:A \rightarrow B$ it is sufficient to find $\forall c,q$ an
injective map $$\alpha_{c,q}: \oplus_{p} A_{c+p}^{p,q}  \otimes 
S^{p,q}V (c+p) 
\rightarrow  \oplus_{\rho} B_{c+\rho}^{\rho,q} \otimes S^{\rho,q}V (c+\rho)
$$ Order $p$ and $\rho $ in decreasing
way, fix a basis of $A_{c+p}^{p,q}$  $\forall p$ 
 and let $\alpha_{c,q}$ be the map s.t.:
{\small  $$M(\alpha_{c,q}) = \left(\begin{array}{c}I \\ 0
\end{array}\right)$$} \hspace{-0.33cm}
Observe that, since $p$ and $\rho$ are ordered in decreasing way and  
$dim(\oplus_{p \geq \tilde{p}} A_{c+p}^{p,q}) \leq dim
(\oplus_{\rho \geq \tilde{p}} B_{c+\rho}^{\rho,q})$, then the entries 
equal to $1$ are ``where $\rho \geq p$''.

{\it i')}  First observe that an $SL(V)$-invariant injective map 
$\alpha: A \rightarrow B$  is s.t.
$0 \rightarrow A \stackrel{\alpha}{\rightarrow} B \rightarrow coker
\, \alpha \rightarrow 0 $ is a minimal free resolution if and only if the maps induced
 by $\alpha$ $$\alpha_{p,q,i} :  A_i^{p,q} \otimes S^{p,q} V (i) \rightarrow
  B_i^{p,q} \otimes S^{p,q} V (i)$$ are zero.

Let $0 \rightarrow A \stackrel{\alpha}{\rightarrow} B \rightarrow coker
\, \alpha \rightarrow 0 $ be a minimal free resolution
 with $ \alpha$ $SL(V)$-invariant injective map.
For any $ c,q, \tilde{p}$, 
let  $$\alpha_{c,q, > \tilde{p}, \geq  \tilde{p}}: \oplus_{p \geq \tilde{p}}
A_{c+p}^{p,q} \otimes S^{p,q}V (c +p)   \rightarrow \oplus_{\rho >
\tilde{p}}  B_{c +\rho}^{\rho,q} \otimes S^{\rho,q}V (c +\rho) $$ the
map induced by $\alpha$.

We can prove that $ M(\alpha_{c,q, > \tilde{p}, \geq  \tilde{p}})$ is injective
as in the implication $\Rightarrow $ of {\it i}, i.e. by considering
$\alpha_{c,q, > \tilde{p} ,\geq  \tilde{p}} \circ \psi_{\tilde{p}}$ 
 with $\psi_{\tilde{p}}$ as above (to see
$\psi_{\tilde{p}} (v \otimes y^{c + \tilde{p}}_{\tilde{p},q})$ 
 is  in $Ker (\alpha)$, use also  that 
the maps  $\alpha_{p,q,i} :  A_i^{p,q} \otimes S^{p,q} V (i) \rightarrow
  B_i^{p,q} \otimes S^{p,q} V (i)$ are zero).
Then \ref{zzz} holds.

The other implication is 
completely analogous to the implication ``if'' of {\it i}.

{\it ii)}
 Suppose 
$\alpha $ is injective.  Fix $q$ and $r$.
Let $\underline{p} = min \{p | \; A^{p,q,r} \neq 0\}$ and
let $$\psi_{ \underline{p} }=  \oplus_{p \geq \underline{p}}  I_{A^{p,q,r}} 
\otimes
\varphi_{p, \underline{p}} : (\oplus_{p \geq \underline{p}}
A^{p,q,r}) \otimes S^{\underline{p},q,r}V (\underline{p}+q+r ) 
\rightarrow \oplus_{p
\geq \underline{p}} A^{p,q,r} \otimes S^{p,q,r}V (p+q+r)$$
It holds  $$\alpha_{q,r}
 \circ \psi_{\underline{p}} =
 (\oplus_{\rho }
I_{B^{\rho,q,r}} \otimes \varphi_{\rho , \underline{p}} )\circ 
(M(\alpha_{q,r}) \otimes I_{S^{\underline{p},q,r}V 
(\underline{p}+q +r)}) $$ 
Let $v \in Ker M(\alpha_{q,r })$. If $v \neq 0$ then 
$\psi_{\underline{p}} (v \otimes 
y^{\underline{p}+q+r}_{\underline{p},q,r})$ 
(see Lemma \ref{xazero} for the definition of 
$y^{\underline{p}+q+r}_{\underline{p},q,r}$) is nonzero (since 
$\psi_{\underline{p}} $ is injective) and it is  in
 $Ker (\alpha)$ (in fact it is in  $Ker(\alpha_{q, r})$
  by the above formula and thus in $Ker (\alpha) $ by the definition of 
$y^{\underline{p}+q+r}_{
\underline{p},q,r}$). Hence we get a contradiction since $\alpha $ is 
injective, therefore $v = 0$. Thus    $ M(\alpha_{q,r})$
is injective; hence   $\alpha_{q,r}$ is injective by Lemma \ref{prelim}.

Now suppose $\alpha_{q,r}$ is injective $\forall q,r$.
Observe  that $\alpha $ is "triangular" with
respect to $q$ and $r$, thus, if 
$\forall q$  the induced map $$\alpha_{q} :\oplus_{p,r}
A^{p,q,r} \otimes S^{p,q,r}V (p+q+r)  \rightarrow \oplus_{\rho,r}
B^{\rho,q,r} \otimes  S^{\rho,q,r}V (\rho+q+r) 
$$ is injective then $\alpha$ is injective; besides $\alpha_q$ is injective
 $\forall q$   if $\forall q,r$  the induced map $\alpha_{q,r}$
is injective.
\hfill  \framebox(7,7)

\medskip

Theorem \ref{fibom} is easily generalizable to ${\bf P}^n$. Obviously
 the statement on minimal resolutions is generalizable to minimal 
free resolutions with two terms of  bundles on ${\bf P}^n$, 
but  for a generic homogeneous bundle on ${\bf P}^n$ with $n \geq 3$
the minimal free resolution has more than two terms.

Finally we observe that 
 Theorem  \ref{fibomintro} and    the following lemma (which will be 
useful also later to study simplicity)
 allow us to parametrize the set of
 homogeneous bundles on ${\bf P}^2$ by a set of sequences
of injective matrices with a certain shape up to the action
 of invertible matrices with a certain shape.

\begin{lemma} \label{soll} i)
Let $E$ and $E'$ be two  homogeneous vector bundles on ${\bf P}^2$ and 
{\small \begin{equation} \label{rismin2}  
a) \;\;\;\; 
 0 \rightarrow R \stackrel{f}{\rightarrow} S \stackrel{g}{\rightarrow}
 E \rightarrow 0 \;\;\;\;\;\;\;\;\;\;\;\; b) \;\;\;\;
 0 \rightarrow R' \stackrel{f'}{\rightarrow} S' \stackrel{g'}{\rightarrow}
 E' \rightarrow 0
\end{equation}} \hspace{-0.3cm} 
be two  minimal  free resolutions. 
 Any map $\eta : E  \rightarrow E'$  induces maps  
$A$ and $B$ s.t. the following  diagram  commutes: 
{\small $$
\begin{array}{ccccccc} 0 \rightarrow  & R &  \stackrel{f}{\rightarrow} &
S
&  \stackrel{g}{\rightarrow}  & E  & \rightarrow 0
\\ & \downarrow A & & \downarrow B & & \downarrow \eta & \\  0 \rightarrow
& R' &  \stackrel{f'}{\rightarrow} & S' &  \stackrel{g'}{\rightarrow} 
 & E'  & \rightarrow 0
\end{array}$$}
\hspace{-0.22cm}
and given $A$ and $B$ s.t. the above diagram commutes we have a map 
 $E \rightarrow E'$. Given $\eta$, 
the maps $A$ and $B$ are unique if and only if $Hom(S,R')=0$.
(In particular 
 if $Hom(S,R')=0$ and we find $A$ and $B$ s.t. the diagram commutes
and $B$ is not a multiple of the identity, we can conclude that $E$ is not 
simple.) 

ii) Let {\small \begin{equation} \label{rismin3} a) \;\;\;\;
0 \rightarrow R \stackrel{f}{\rightarrow} S \stackrel{g}{\rightarrow}
 E \rightarrow 0 \;\;\;\;\;\;\;\;\;\;\;\;\;\; b) \;\;\;\;
 0 \rightarrow R \stackrel{f'}{\rightarrow} S \stackrel{g'}{\rightarrow}
 E \rightarrow 0 \end{equation} } \hspace{-0.3cm}  be 
 two  minimal  free resolutions with $SL(V)$-invariant maps
 of a homogeneous bundle  $E$ on ${\bf P}^2 = {\bf P}(V) $;
 then there exist  $SL(V)$-invariant 
automorphisms  $A : R \rightarrow R$, $B : S \rightarrow S$
  s.t. the following diagram commutes:  {\small $$
\begin{array}{ccccccc} 0 \rightarrow  & R &  \stackrel{f}{\rightarrow} &
S &  \stackrel{g}{\rightarrow}  & E  & \rightarrow 0
\\ & \downarrow A & & \downarrow B & & \downarrow Id & \\  0 \rightarrow
& R &  \stackrel{f'}{\rightarrow} & S &  \stackrel{g'}{\rightarrow} 
 & E  & \rightarrow 0
\end{array}$$} 
\end{lemma}

{\it Proof.} {\it i)}  The composition  $S \stackrel{g}{\rightarrow} 
 E \stackrel{\eta}{\rightarrow} E'$ 
can be lifted to a  map $B: S \rightarrow S'$
by the exact sequence 
$ 0 \rightarrow Hom (S,R')  \rightarrow Hom (S,S')  \rightarrow Hom 
(S,E') \rightarrow 0$ (obtained  by applying  $Hom (S, \cdot) $ to 
(\ref{rismin2})b).
By the exact  sequence
$Hom (R,R') \rightarrow Hom (R,S') \rightarrow Hom(R,E')$ (obtained by  
applying $Hom (R, \cdot) $  to (\ref{rismin2})b)
    the composition $R \stackrel{f}{\rightarrow}
 S \stackrel{B}{\rightarrow} S' $, which goes to zero in $Hom(R,E')$,
can be lifted to a  map  $A: R \rightarrow R'$.   

{\it ii)} As in {\it i}
 we can prove there exist $SL(V)$-invariant maps $ B : S 
\rightarrow  S $ and  $ B' : S \rightarrow  S $ s.t. the following diagrams 
commute
{\small $$
\begin{array}{ccccccc} 0 \rightarrow  & R &  \stackrel{f}{\rightarrow} &
S &  \stackrel{g}{\rightarrow}  & E  & \rightarrow 0
\\ &  & & \downarrow B & & \downarrow Id & \\  0 \rightarrow
& R &  \stackrel{f'}{\rightarrow} & S &  \stackrel{g'}{\rightarrow} 
 & E  & \rightarrow 0
\end{array} \;\;\;\;\;\;\;\;\;\;\;\;
\begin{array}{ccccccc} 0 \rightarrow  & R &  \stackrel{f'}{\rightarrow} &
S &  \stackrel{g'}{\rightarrow}  & E  & \rightarrow 0
\\ &  & & \downarrow B' & & \downarrow Id & \\  0 \rightarrow
& R &  \stackrel{f}{\rightarrow} & S &  \stackrel{g}{\rightarrow} 
 & E  & \rightarrow 0
\end{array}$$}  
\hspace{-0.22cm} 
 then we get the commutative diagram {\small 
$$\begin{array}{ccccccc} 0 \rightarrow  & R &  \stackrel{f}{\rightarrow} &
S &  \stackrel{g}{\rightarrow}  & E  & \rightarrow 0
\\ &  & & \;\;\;\;\downarrow B' \circ B & & \downarrow Id & \\  0 \rightarrow
& R &  \stackrel{f}{\rightarrow} & S &  \stackrel{g}{\rightarrow} 
 & E  & \rightarrow 0
\end{array}$$}
\hspace{-0.3cm} Since also $Id:S \rightarrow S$  lifts  
$Id : E \rightarrow E$ in the above diagram (i.e. $Id \circ g = g \circ Id$), 
  there exists a map $a : S \rightarrow R$ 
s.t. $ B' \circ B - Id = f \circ a $; hence, by the minimality of 
(\ref{rismin3}){\it a},
  $det (B' \circ B ) = 1 + P$ where $P$ is a polynomial 
without terms of degree $0$; since $det (B' \circ B) : det (S) \rightarrow 
det (S)  $ must be homogeneous, $P=0$; hence  $B' \circ B$ is invertible
and thus $B$ is invertible. 
\hfill  \framebox(7,7)

\section{Quivers}

We recall now the main definitions and results on  
 quivers and representations of  quivers associated to
 homogeneous bundles introduced by Bondal and Kapranov in \cite{B-K}.
The quivers will allow
 us to handle well and ``to make explicit'' the homogeneous 
subbundles of a homogeneous bundle.


\begin{defin} (See  \cite{Sim}, \cite{King}, \cite{Hil1}, \cite{G-R}.)
A {\bf quiver} is an oriented  graph ${\cal Q}$ with the set ${\cal Q}_0$ of 
vertices (or points) and the set  ${\cal Q}_1$ of arrows.

A {\bf path} in ${\cal Q}$ is  a formal composition of arrows $ \beta_m ...
\beta_1$ where the source of an arrow $\beta_i $
is the sink of the previous arrow $\beta_{i-1}$.

A {\bf relation} in ${\cal Q}$ is a linear form $\lambda_1 c_1+...+
\lambda_r c_r$ where $c_i$ are paths in ${\cal Q}$ with a common source and
 a common sink and $\lambda_i \in {\bf C}$.

A {\bf representation of a quiver} ${\cal Q} =({\cal Q}_0, 
{\cal Q}_1)$, or {\bf ${\cal Q}$-representation},
  is the couple of a set of vector 
spaces $\{X_i\}_{i \in {\cal Q}_0} $ and of a set of linear maps 
$ \{\varphi_{\beta} \}_{\beta \in {\cal Q}_1}$ where  $\varphi_{\beta} : X_i 
\rightarrow X_j$ if $\beta$ is an arrow from $i$ to $j$.

A {\bf representation of a  quiver ${\cal Q}$ with relations ${\cal R}$} 
is a ${\cal Q}$-representation s.t. 
$$\sum_j \lambda_j \varphi_{{\beta}^j_{m_j}} ...\varphi_{{\beta}^j_{1}}=0$$ 
for every $ \sum_j \lambda_j \beta^j_{m_j} ...\beta^j_{1} \in {\cal R}$.

Let $(X_i, \varphi_{\beta})_{i \in {\cal Q}_0,\; \beta \in {\cal Q}_1}$ and 
$(Y_i, \psi_{\beta})_{i \in {\cal Q}_0,\; \beta \in {\cal Q}_1}$ be two 
representations of a quiver ${\cal Q}= ({\cal Q}_0 , {\cal Q}_1)$.
A {\bf morphism} $f $ from $(X_i, \varphi_{\beta})_{i \in {\cal Q}_0,\; 
\beta \in {\cal Q}_1}$ to
$(Y_i, \psi_{\beta})_{i \in {\cal Q}_0,\; \beta \in {\cal Q}_1}$ is a
 set of linear maps $f_i : X_i \rightarrow Y_i $, $i \in {\cal Q}_0$ s.t.,
 for
 every  $\beta \in {\cal Q}_1$,  $\beta$  arrow from $ i$ to 
 $j$, the following diagram is commutative: 
{\small $$\begin{array}{ccc} X_i   
 & \stackrel{f_i}{\longrightarrow}
 &  Y_i \\ 
\varphi_{\beta} \downarrow & &  \downarrow \psi_{\beta}
\\   X_j 
 &  \stackrel{f_j}{\longrightarrow} &  Y_j
 \end{array}$$ } \hspace{-0.3cm}
A morphism $f$ is injective if the $f_i$ are injective.
\end{defin}

\begin{notat}
We will say that a representation $(X_i, \varphi_{\beta})_{i \in 
{\cal Q}_0,\; \beta \in {\cal Q}_1}$ 
of a quiver ${\cal Q}= ({\cal Q}_0 , {\cal Q}_1 )$ has 
{\bf multiplicty}  $m$  in a point $i$  of ${\cal Q}$  if $dim X_i =m$.
 
The  {\bf support} (with multiplicities)
 of a representation of a quiver ${\cal Q}$
is the subgraph
of ${\cal Q}$  constituted by the points  of multiplicity  $\geq 1$ and the
 nonzero arrows (with the multiplicities
 associated to every point of the subgraph).
\end{notat}

We recall now from  \cite{B-K}, \cite{Hil1}, \cite{Hil2}
the definition of a quiver ${\cal Q}$ 
 s.t.  the category of the homogeneous bundles on ${\bf P}^2$ 
is equivalent
to the category of finite dimensional representations of ${\cal Q}$
with some relations ${\cal R}$. Bondal and Kapranov defined  such a quiver
in a more general setting but we recall
such a  construction only for ${\bf P}^2$. See also \cite{O-R}.


\medskip

First some notation.
Let $P$ and $R$ be the following  subgroup of $SL(3)$:
$$ P= \{ {\scriptsize
 \left(  \begin{array}{ccc} a & b & c \\ 0 & d & e \\ 0 & f & g
 \end{array}     \right)} \in SL(3)\} \;\;\;\;\;\;\;\;\;\;\;\;\;
 R= \{ {\scriptsize
\left(  \begin{array}{ccc} a & 0 & 0 \\ 0 & d & e \\ 0 & f & g
 \end{array}     \right)} \in SL(3)\}$$
Observe that $R$ is reductive.
We can see ${\bf P}^2$ as $$ {\bf P}^2 = SL(3) /P$$ 
($P$ is the stabilizer of   $[ 1 : 0 : 0]$). 
Let $p$ and $r$ be the Lie algebras
 associated respectively to $P$ and $R$. 
 Let $n$ be the Lie algebra 
$$n= \{ {\scriptsize
\left(  \begin{array}{ccc} 0 & x & y \\ 0 & 0 & 0 \\ 0 & 0 & 0
 \end{array}     \right)} | \; x,y \in {\bf C}\}$$ Thus $ p=r \oplus n$
(Levi decomposition).

We recall  that  the homogeneous bundles on  $ {\bf P}^2 = SL(3) /P$ are 
given by the representations of $P$ (this bijection is given by taking the 
fibre over  $[ 1 : 0 : 0]$  of a homogeneous vector bundle); by
 composing the projection from $P$ to $R$ 
with a representation of $R$ we get a representation of $P$ and the set of
the homogenous bundles obtained in this way from the irreducible 
representations of $R$  are $$\{S^l Q (t) 
| \; l \in {\bf N}, \; t \in {\bf Z}\} $$ 
where $Q = T_{{\bf P}^2 }(-1)$;

\medskip

\begin{defin}
From now on ${\cal Q}$ will be the following quiver:

$\bullet $ let $${\cal Q}_0 = \{irreducible \; representations\; of\; R \}=
\{S^l Q (t) 
| \; l \in {\bf N}, \; t \in {\bf Z}\}= $$
$$ = \{dominant \; weights \; of \; r\}$$
$\bullet $
let ${\cal Q}_1$ be defined in the following way: there is an arrow from 
$\lambda $ to $\mu$, $\lambda, \mu \in {\cal Q}_0$,  if and only if  $ n \otimes 
\Sigma_{\lambda} \supset \Sigma_{\mu}$, where $\Sigma_{\lambda}$ denotes 
the representation  of $r$ with dominant weight $\lambda$.
\end{defin}

\begin{lemma}  The adjoint representation of $p$ on $n$ 
 corresponds to  $Q(-2)= \Omega^1$, 
 more precisely: 
let $\rho : p \rightarrow gl(n)$ be the following  representation: 
$$\rho (B ) h( {\small \left( \hspace{-0.1cm}
 \begin{array}{c} x \\ y \end{array}   \hspace{-0.1cm}   \right)})
=  Bh  {\small \left(   \hspace{-0.1cm}
\begin{array}{c} x \\ y \end{array}    \hspace{-0.1cm}  \right)}
 -h  {\small \left(  \hspace{-0.1cm} \begin{array}{c} x \\ y \end{array} 
\hspace{-0.1cm}    \right) } B$$
where  $h :   {\bf C}^2 \rightarrow n $ is the isomorphism
$ {\small \left(  \hspace{-0.1cm} 
 \begin{array}{c} x \\ y \end{array}    \hspace{-0.1cm} \right)}  \mapsto
{\scriptsize  \left(  \begin{array}{ccc} 0 & x & y \\ 0 & 0 & 0 \\ 0 & 0 & 0
 \end{array}  \right)}  $; $Q(-2)$ is 
the homogeneous bundle whose fibre as $p$-representation is $ n$.
\end{lemma}

{\it Proof.}
Observe that if  $B =  {\scriptsize 
\left(  \begin{array}{ccc} a & b & c \\ 0 & d & e  \\ 0 & f  & g
 \end{array}     \right)} \in p $ and
   $ A =  {\scriptsize \left(  \begin{array}{cc}  d & e \\ f & g
 \end{array}     \right)}$ then 
$$h^{-1} (Bh {\small \left( \hspace{-0.1cm}
  \begin{array}{c} x \\ y \end{array}   \hspace{-0.1cm}  \right)}
 -h  {\small \left( \hspace{-0.1cm} \begin{array}{c} x \\ y \end{array}   
\hspace{-0.1cm}   \right)} 
B) = (a \; Id - {}^t A)  {\small \left(  \hspace{-0.1cm}
\begin{array}{c} x \\ y \end{array}  \hspace{-0.1cm}
   \right)}; $$  
since the representation $p \rightarrow gl(1) $,
 $B \mapsto a \; Id$ corresponds to the bundle $ {\cal O} (-1)$ and the  
representation $p \rightarrow gl(2) $,
 $B \mapsto - {}^t A$ corresponds to the bundle $ Q^{\vee}= Q(-1)$,
 we conclude. \hfill  \framebox(7,7)

\medskip

By the previous remark,
 if $\Sigma_{\lambda}$ is the representation corresponding to
$ S^l Q (t) $ and $\Sigma_{\mu}$ is the representation corresponding to
$ S^{l'} Q (t') $,   the condition  $ n \otimes  
\Sigma_{\lambda}   \supset \Sigma_{\mu}$ is equivalent to the fact 
$S^{l'} Q (t')$
is a direct summand of  $ Q (-2) \otimes S^l Q (t) $ and this is true
if and only if $ (l', t') =(l-1, t-1)$ or  $ (l', t') =(l+1, t-2)$
 (we recall that, by the Euler sequence, $\wedge^2 Q = 
{\cal O}(1)$).

Thus our quiver has three connected components ${\cal Q}^{(1)}$,  
${\cal Q}^{(2)}$,  ${\cal Q}^{(3)}$  
(given by the congruence class modulo $3/2$ of the slope of the homogeneous
 vector bundles corresponding to the points of the connected component); 
the figure shows one of them (the one whose points correspond to the 
bundles with  $\mu \equiv 0 \; \,mod \,\;(3/2)$):
we identify the points of every connected component ${\cal Q}^{(j)}$  
of ${\cal Q}$    with a subset of ${\bf Z}^2$   for convenience.

\begin{center}
\includegraphics[scale=0.4]{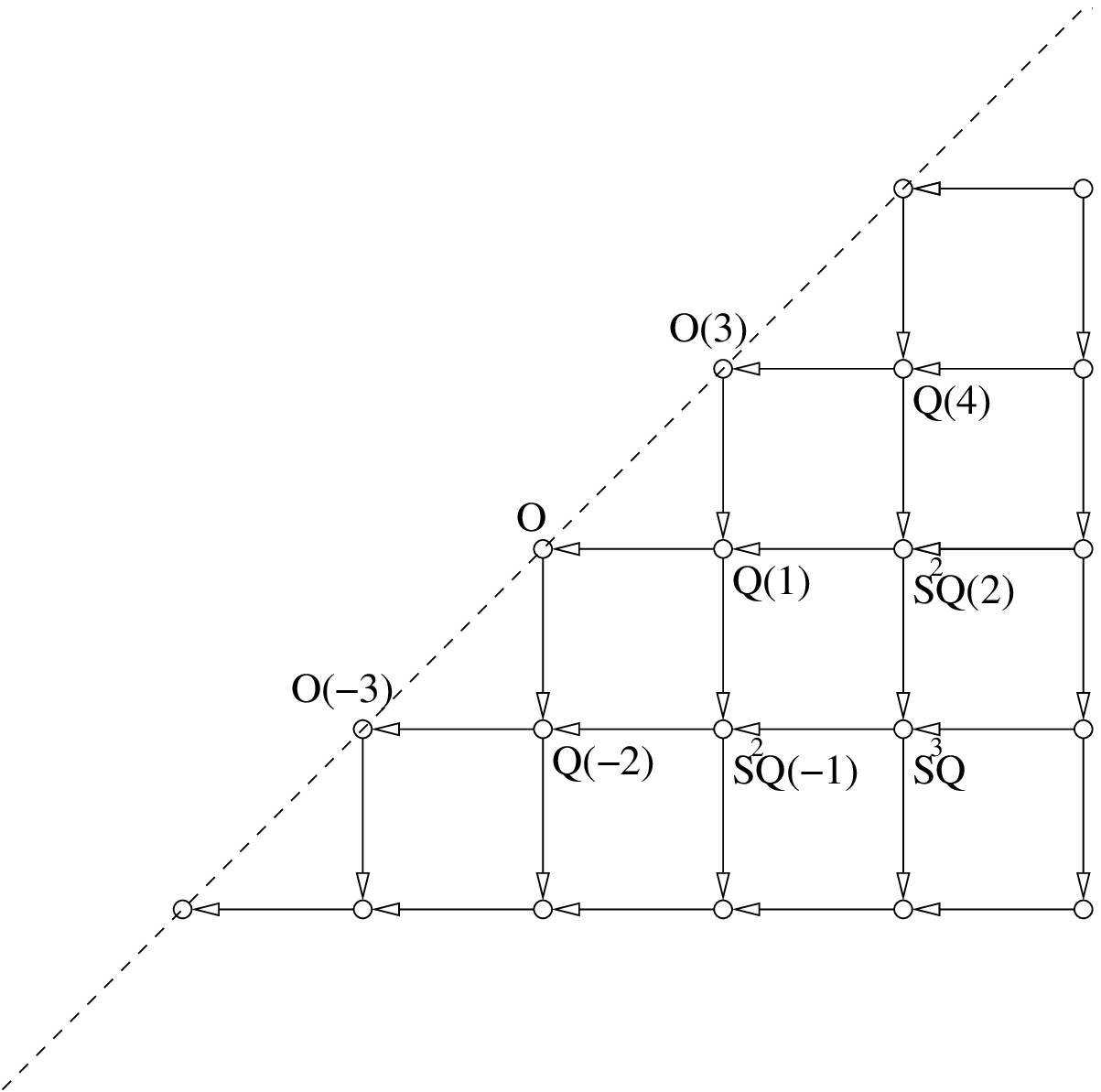}

\end{center}

\smallskip

\begin{defin} \label{relazioni}
Let ${\cal R}$ be the set of relations on ${\cal Q}$  given
 by the commutativity of the squares, i.e. (denoting  by $\beta_{w,v}$ the 
arrow from $v$ to $w$):
 $$ \beta_{ (x -1, y -1), (x -1 , y) }
 \beta_{  (x -1, y),  (x , y) } - 
\beta_{ (x -1 , y -1),  (x , y -1)  }
\beta_{ (x , y -1), (x , y)   }$$
$\forall (x,y)   \in 
{\cal Q}^{(j)} \subset {\bf Z}^2 $
 for some $j$, s.t. $ (x -1, y ) \in {\cal Q}$ and  $$
\beta_{ (x -1 , y -1),  (x , y -1)  }
\beta_{ (x , y -1), (x , y)   }$$ $\forall (x,y)   \in 
{\cal Q}^{(j)}  \subset {\bf Z}^2 $ for some $j$, s.t.  
$ (x -1, y ) \not\in  {\cal Q}$.
\end{defin}

\begin{center}
\includegraphics[scale=0.4]{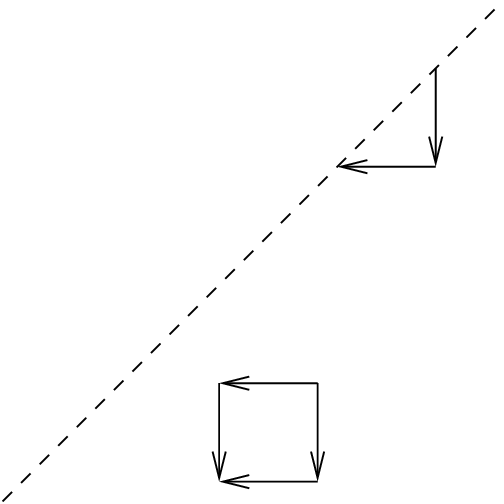}
\end{center}

\vspace{0.05cm}

\begin{defin} \label{rapp}
Let $E$ be a homogeneous vector bundle on ${\bf P}^2$. The {\bf 
 ${\cal Q}$-representation  associated to $E$} is the following
 (see \cite{B-K}, \cite{Hil1}, \cite{Hil2}, \cite{O-R}): 

consider the fiber $E_{[1:0:0]}$  of $E$ on $[1:0:0]$; 
it is a representation  of $p$; as $r$-representation we have 
 $$E_{[1:0:0]}  =
\oplus_{\lambda \in {\cal Q}_0} X_{\lambda} \otimes  \Sigma_{\lambda} $$ 
for some vector spaces $X_{\lambda}$;
 we associate to $ \lambda \in {\cal Q}_0$ the vector space $X_{\lambda}$;
we fix  $\forall \lambda $ a dominant  weight vector $v_{\lambda} \in
 \Sigma_{\lambda}$ and $\eta_1, \eta_2$ eigenvectors of the 
$p$-representation $n$; let  $\psi_1, \psi_2$ be their weights  respectively; 
let $i$  be  s.t. $\lambda+\psi_i = \mu$;
we associate  to an arrow $\lambda \rightarrow \mu$ a map $f: X_{\lambda}
\rightarrow X_{\mu}$ defined in the following way: consider the composition
$$
 \Sigma_{\lambda} \otimes n
 \otimes X_{\lambda} \longrightarrow \Sigma_{\mu} \otimes X_{\mu}$$

 given by the action of $n$ over $E_{[1:0:0]}$
 followed by projection; 
it maps $ v_{\lambda} \otimes 
\eta_i \otimes v $ to $v_{\mu} \otimes w$; we define $f(v)=w$
(it does not depend on the choice of the dominant weight vector).
\end{defin}

\begin{theorem} \label{BKH}
 {\bf  (Bondal, Kapranov, Hille)}  \cite{B-K}, \cite{Hil1}, 
\cite{Hil2}, \cite{O-R}.  The
 category of the homogeneous bundles on ${\bf P}^2$ is equivalent
to the category of finite dimensional representation of the quiver ${\cal Q}$
with the relations ${\cal R}$.
\end{theorem}

\vspace{0.05cm}

Observe that in Def. \ref{rapp}
with respect to Bondal-Kapranov-Hille's
 convention in  \cite{B-K}, \cite{Hil1}, 
\cite{Hil2}, we preferred to invert the arrows 
in order that  an injective $SL(V)$-equivariant  map 
of bundles   corresponds to an injective morphism of 
${\cal Q}$-representations.
 For example ${\cal O}  $ injects in $ V (1) $ whose support 
is the arrow  from  $Q(1) $ to ${\cal O}$. 
$$\begin{array}{c} \circ \hspace{-0,15cm}  \longleftarrow  \hspace{-0,2cm}
 \circ \\ \hspace*{0,5cm}
{\cal O} \hspace{0,4cm} Q (1) 
 \end{array} $$

\vspace{0.05cm}

\begin{notat} $\bullet $ 
We will often speak of the {\bf ${\cal Q}$-support} 
of a homogeneous bundle $E$ 
instead of the  support with multiplicities 
of the ${\cal Q}$-representation of $E$ and we will 
denote it by ${\cal Q}$-$supp(E)$.

 $\bullet $
The word ``{\bf rectangle}'' will denote the subgraph with multiplicities 
of ${\cal Q}$ given by the subgraph of ${\cal Q}$  included
 in a rectangle whose sides are unions of  arrows of ${\cal Q}$, with the
multiplicities  of all its points equal to $1$.

The word ``{\bf segment}'' will denote a rectangle with base or height equal
 to $0$.

 $\bullet $ If $A$ and $B$ are two subgraphs of ${\cal Q} $, $A \cap B$ is the
 subgraph of ${\cal Q}$ whose vertices and arrows are the vertices and arrows
both of $A$ and of $B$; $A -B$ is the subgraph of ${\cal Q} $ whose vertices 
are the vertices of $A$ not in $B$ 
   and the arrows are the arrows of $A$ joining two vertices of $A -B$.

\end{notat}

\begin{rem} \label{Spq}  (\cite{B-K}) {\rm 
The ${\cal Q}$-support of $S^{p,q}V$ is  
a rectangle  as in the figure: 
\begin{center}
\includegraphics[scale=0.33]{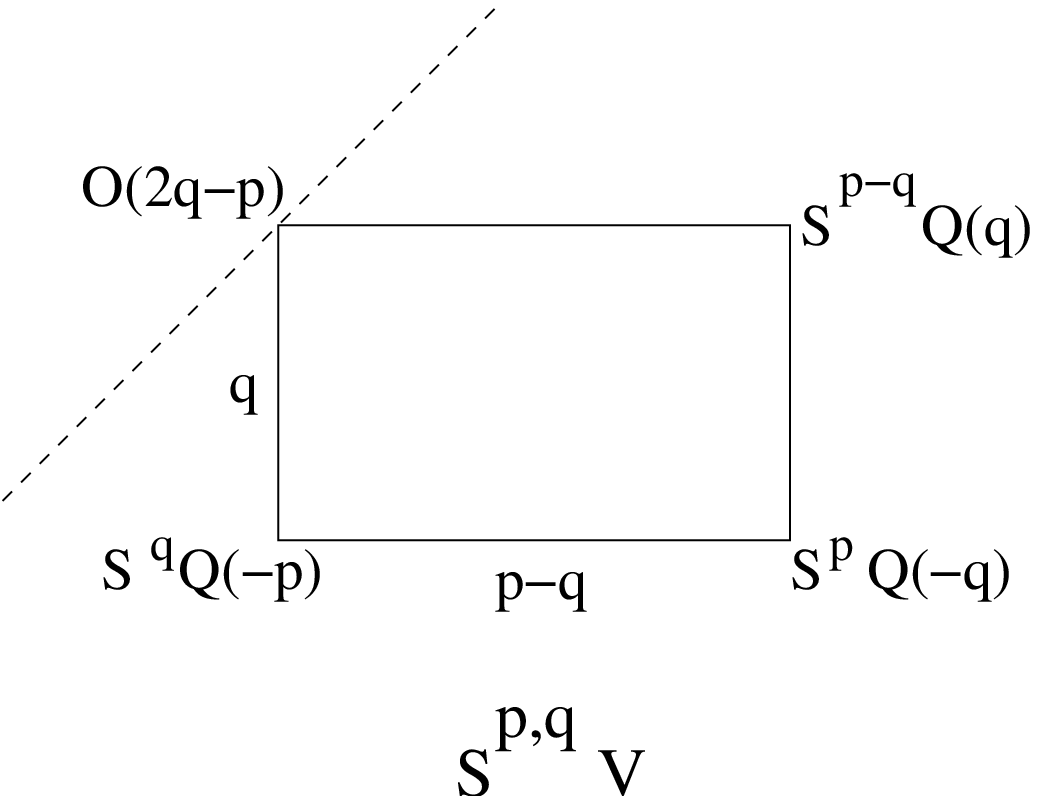}
\end{center}
In fact:
  by  the Euler sequence $ S^{p,q}V
= S^{p,q} ({\cal O} (-1) \oplus Q )  $ as $R$-representation; by  the
 formula of a Schur functor applied to a direct sum (see
\cite{F-H}, Exercise 6.11) we get 
$$S^{p,q} V = \oplus S^{\lambda} Q \otimes S^m {\cal O}
(-1)   $$ as $R$-representations, 
 where the sum is performed on $ m \in {\bf N}$ and on $ \lambda $ Young 
diagram obtained  from the Young diagram of $(p,q)$ by taking off $m$ boxes not
 two in the same column; thus  $$
S^{p,q} V  =\oplus_{0 \leq m_1 \leq p-q , \;\;
0 \leq m_2 \leq q   } S^{p-m_1, q-m_2} Q (-m_1 -m_2)$$ 
Finally to show the maps associated to the arrows in the rectangle are nonzero 
we can consider on the set of the vertices of the rectangle the following
 equivalence relation: $P \sim Q$ if and only if there exist two paths
 with $P$ and $Q$ respectively as sources and 
common sink s.t. the  map associated 
to any arrow of the two paths is nonzero; if a map associated to an arrow 
(say from $P_1 $ to $P_2$)   of  the rectangle    is zero then, by the 
``commutativity of the squares'', precisely by the relations in Definition 
\ref{relazioni},
 there would be at least two equivalence classes (the class of $P_1$ 
and the class  of $P_2 $), but this is impossible  by the irreducibility 
of $S^{p,q} V$. }
\end{rem}

\section{Some lemmas and notation}

In this section we 
study equivariant maps between some homogeneous bundles on
 $ {\bf P}^2$ by using the language of the quivers  and
we collect some technical notation and lemmas, which will be useful in the 
next section to study homogeneous subbundles (in particular their slope) 
of homogeneous bundles and then to study stability.

\subsection{${\cal Q}$-representation of kernels and images}

\begin{lemma} \label{palle} 
Let  $\varphi: S^{p,q} V \rightarrow S^{p+s_1, q+s_2 , s_3}V 
(s_1+s_2 +s_3) $ be  an 
$SL(V)$-invariant nonzero map;
then 

i) ${\cal Q}$-$supp(Ker(\varphi)) = 
{\cal Q}$-$supp(S^{p,q} V) - {\cal Q}$-$supp( S^{p+s_1, q+s_2 , s_3}V 
(s_1+s_2 +s_3)) $

ii) ${\cal Q}$-$supp(Im(\varphi)) = 
{\cal Q}$-$supp(S^{p,q} V) \cap {\cal Q}$-$supp( S^{p+s_1, q+s_2 , s_3}V 
(s_1+s_2 +s_3)) $

iii) ${\cal Q}$-$supp(Coker(\varphi)) = 
{\cal Q}$-$supp( S^{p+s_1, q+s_2 , s_3}V 
(s_1+s_2 +s_3)) - {\cal Q}$-$supp(S^{p,q} V) $.
\end{lemma}

{\it Proof.} Consider the $P$-invariant (and thus $R$-invariant) map 
$\varphi_{[1:0:0]}$
  induced by $\varphi$ on the fibers on $[1:0:0]$; it is the 
morphism from the ${\cal Q}$-representation of $S^{p,q} V $ to the 
${\cal Q}$-representation of  $ S^{p+s_1, q+s_2 , s_3}V 
(s_1+s_2 +s_3) $.

Obviously the $R$-representations corresponding to the vertices 
in the  ${\cal Q}$-support of $S^{p,q} V$ and not in the ${\cal Q}$-support of 
$ S^{p+s_1, q+s_2 , s_3}V  (s_1+s_2 +s_3) $ are in 
$Ker \varphi_{[1:0:0]}$.

We have to show that 
 the $R$-representations corresponding to the vertices 
both in the  ${\cal Q}$-support of $S^{p,q} V$ and  in the 
${\cal Q}$-support of 
$ S^{p+s_1, q+s_2 , s_3}V  (s_1+s_2 +s_3) $ are not  in 
$Ker \varphi_{[1:0:0]}$. We call $I$ the set of such vertices. 
If the $R$-representation
 corresponding to an element of $I$  is in  $Ker \varphi_{[1:0:0]}$, 
then also the $R$-representation corresponding to another element of $I$  
 is in  $Ker \varphi_{[1:0:0]}$: in fact 
by the commutativity of the diagram 
in the definition of morphism of representations of a quiver, if 
the $R$-representation
 corresponding to  $ \lambda \in I$  is in  $Ker \varphi_{[1:0:0]}$,
 then also the $R$-representation 
corresponding to any element of $I$ linked to $\lambda $ by an arrow  
 is in  $Ker \varphi_{[1:0:0]}$ and we conclude  since $\forall \lambda_1 ,
\lambda_2 \in  I$ there is a path of the quiver joining $\lambda_1 $ and
$\lambda_2$.
 Thus either any $R$-representation 
corresponding to an element of $I$ is in  $Ker \varphi_{[1:0:0]}$ or any 
$R$-representation  corresponding to an element of $I$ is not 
 in  $Ker \varphi_{[1:0:0]}$. But the last case is impossible because 
$\varphi$ is nonzero. 
 \hfill  \framebox(7,7)

\begin{cor} \label{palle2}
Let  $\varphi: 
S^{p,q} V \rightarrow \oplus_{s_1,s_2,s_3} S^{p+s_1, q+s_2 , s_3}V 
(s_1+s_2 +s_3) $ (where the sum is on a finite subset of ${\bf N}^3$)
be an $SL(V)$-invariant nonzero map. 

Then  ${\cal Q}$-$supp(Ker(\varphi))= 
{\cal Q}$-$supp(S^{p,q} V) - {\cal Q}$-$supp( \oplus 
S^{p+s_1, q+s_2 , s_3}V (s_1+s_2 +s_3)) $.
\end{cor}

\begin{lemma} \label{4ter} {\bf (Four Terms Lemma.)}
 On ${\bf P}^2 ={\bf P}(V) $
we have the following exact sequence: $$0 \rightarrow S^{q+s -1,q}
V (-p +q -1 +s) \rightarrow S^{p,q} V \rightarrow S^{p,q+s} V (s)
\rightarrow S^{p-q-1,s-1} V (q+1+s) \rightarrow 0 $$ where the
maps are   $SL(V)$-invariant nonzero
maps (they are unique up to multiples).
 (Observe that $ S^{p-q-1,s-1} V (q+1+s) \simeq S^{p,q+s, q+1} V (q+1+s)$.) 
\end{lemma}

{\it Proof.} In the figure we  show the sides of
the ${\cal Q}$-supports  of $S^{p,q}V$ and $S^{p,q+s}V (s)$:

\begin{center}
\includegraphics[scale=0.37]{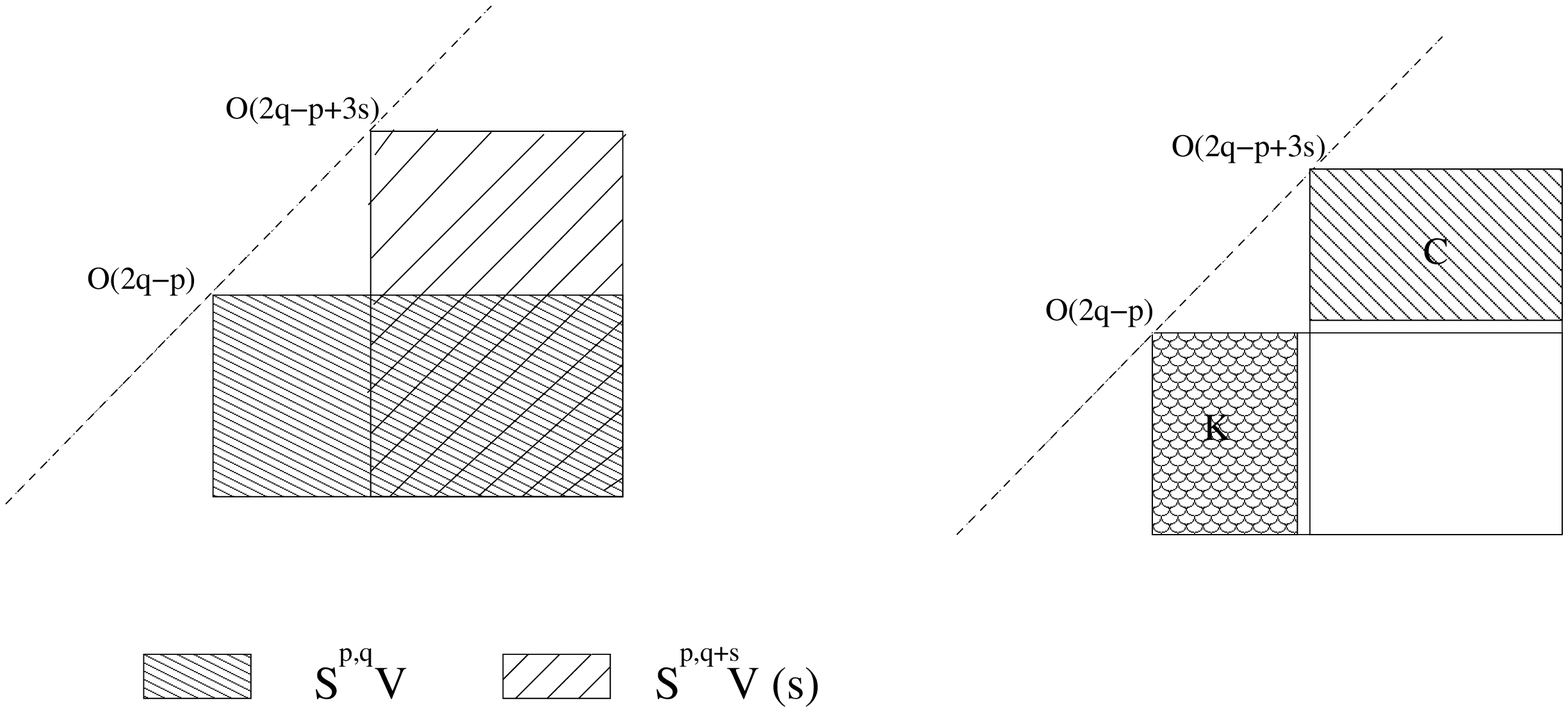}
\end{center}

these supports  are rectangles 
(see Remark \ref{Spq}); thus by Lemma \ref{palle} 
the ${\cal Q}$-support 
of the kernel of   $S^{p,q} V \rightarrow
S^{p,q+s} V (s)$ is the rectangle $K$, which is the 
${\cal Q}$-support of  $ S^{q+s -1,q} V (-p +q -1 +s)$
and the  ${\cal Q}$-support 
of the cokernel of   $S^{p,q} V \rightarrow
S^{p,q+s} V (s)$ is the rectangle $C$, which is the 
support of $S^{p-q-1,s-1} V (q+1+s)$. 
\hfill  \framebox(7,7)

\subsection{Some calculations on the slope}

\begin{rem}  \label{gio} {\rm i) The first Chern class
  of a homogeneous bundle $E$ can be calculated as the 
sum of the first Chern classes of the irreducible bundles corresponding to 
the vertices of the ${\cal Q}$-support of $E$ multiplied 
by the multiplicities. The rank of $E$ is  the 
sum of the ranks of the irreducible bundles corresponding to 
such vertices  multiplied by the multiplicities.

We will often speak of the slope
 (resp. first Chern class, rank) of  a 
graph with multiplicities  instead of
 the slope (resp. first Chern class, rank) of the  vector 
bundle whose ${\cal Q}$-support is  that graph with multiplicities.

ii) Suppose  the set of the vertices of the 
${\cal Q}$-support of $E$ is the disjoint union of the 
vertices of the supports of two ${\cal Q}$-representations $A$ and $B$;
if $\mu(A)= \mu(B) $ then $\mu(E)= \mu(A) = \mu(B)$, if  
 $\mu(A) < \mu(B) $ then $\mu(A) < \mu(E) < \mu(B)$.

iii) We recall that the rank of $S^l Q(t) $ is $l+1$ and its first Chern class
is $ (l+1) (l/2 +t)$. }
\end{rem}

\begin{lemma} \label{murettangolo} 
Let $R$ be a rectangle of base $h$, height $k$ and $S^l Q (t)$ as  the highest
 vertex of the  left side. Then 
 $$\mu (R)= \frac{(h+1)(k+1)[\frac{h^2 -k^2}{2} +h (l+\frac{t}{2} +1) +
k(\frac{t}{2} - \frac{l}{2} -1)+ 
(l+1)(\frac{l}{2} +t)) ] }{ (h+1) (k+1)  (l + \frac{h+k}{2} +1)}  $$
(where the numerator is the first Chern class and the denominator is the rank.)
\end{lemma}

{\it Proof.} Left to the reader. 
\hfill  \framebox(7,7)

\begin{lemma} \label{segmenti} Let $S$ be 
 a  horizontal (resp. vertical) segment in ${\cal Q}^{(j)}$
for some $j$ and let $S'$ be   obtained translating 
$S$ by $(0,1)$ (resp $(1,0)$) in  ${\cal Q}^{(j)} \subset {\bf Z}^2$. 
Then $\mu(S) < \mu(S')$.
\end{lemma}

\begin{center}
\includegraphics[scale=0.28]{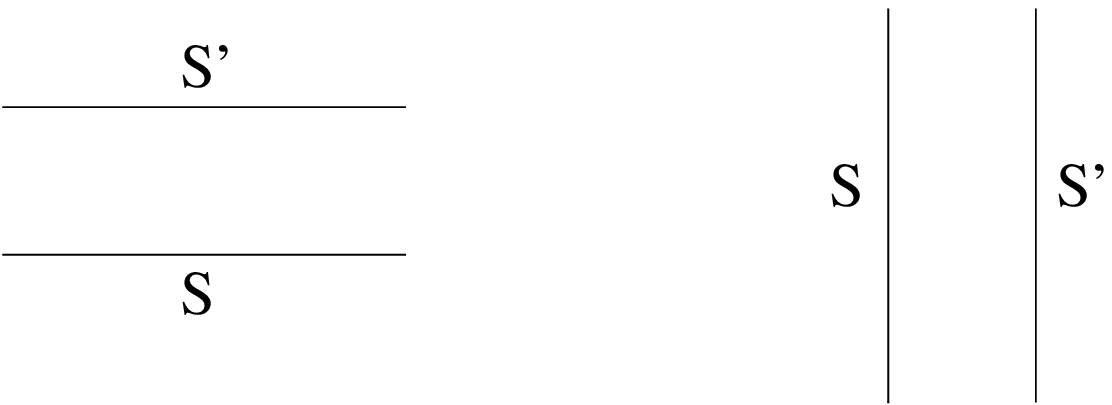}
\end{center}

{\it Proof.}  Suppose $S$ is horizontal, of length $h$
and  its  first vertex from left is $S^{l} Q(t) $. By Lemma 
\ref{murettangolo}, 
{\small
$$ \mu (S') =\frac{h^2 + 2 h (l+1) +l(l+1) }{2 l +2 +h} +t$$ 
$$ \mu (S) =\frac{h^2 + 2 h (l+2) +(l+2)(l+1) }{2 (l+1) +2 +h} +t -2$$} 
\hspace{-0.26cm} 
and $ \mu (S') > \mu (S)$  is easy to check. The case of vertical segments
 is similiar. 
\hfill  \framebox(7,7)

\begin{lemma} \label{rettangoli}
  i) Let $U$ and $U'$ be two rectangles with the same base 
with the vertical sides lined up and $U$ above $U'$.  Then 
$\mu (U') < \mu (U)$.  

ii) Let $W$ and $W'$ be two rectangles with the same height  
and  the horizontal sides lined up and $W$ at the right side of  $W'$.  Then 
$\mu (W') < \mu (W)$.

iii) Let $R$ be a  rectangle and $ R'$ a subrectangle of $R$ with the same
 base and with the lower side equal to the lower side of $R$. Then   $\mu(R') 
< \mu (R) $. 

iv) Let $T$ be a  rectangle and $ T'$ a subrectangle of $T$ with the same
 height and with the left side equal to the left side of $T$. Then   $\mu(T') 
< \mu (T) $.

\begin{center}
\includegraphics[scale=0.3]{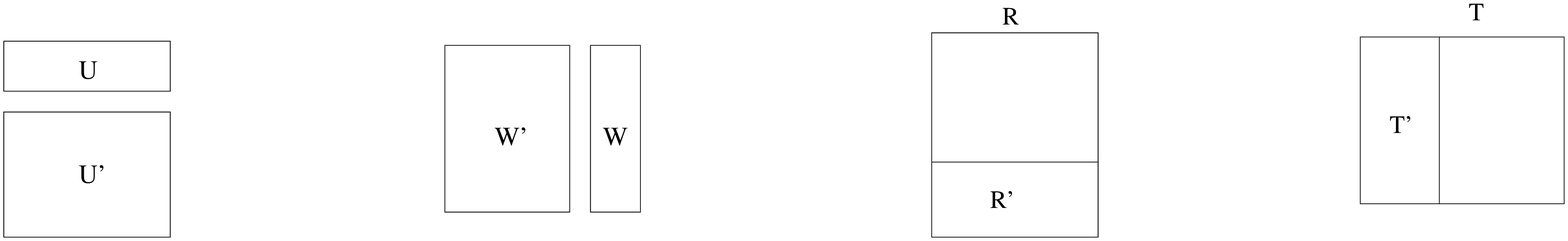}
\end{center}
\end{lemma}

{\it Proof.} {\it i)}
 By Lemma \ref{segmenti} and Remark \ref{gio} 
the slope of a rectangle  is between the slope 
of the lower side and the slope of the higher side, thus again by Lemma 
\ref{segmenti} we conclude.

{\it ii)} Analogous to {\it i}.

{\it iii)} and {\it iv)} follow from {\it i} and {\it ii}. 
\hfill  \framebox(7,7)

\subsection{Staircases}

In this subsection we introduce particular ${\cal Q}$-representations, called
 ``staircases''. Their importance is due to the fact that they are the 
${\cal Q}$-supports of  the homogeneous subbundles 
 of the homogeneous bundles whose
${\cal Q}$-supports are rectangles (in particular  of the  
trivial homogeneous bundles).

\begin{rem} \label{subquiv} {\rm 
Let $E$ be a homogeneous bundle on ${\bf P}^2$ and $F$ be a homogeneous
 subbundle. 
Let   $S$ and $S'$  be the ${\cal Q}$-supports
of $E$ and   $F$ respectively.  By Theorem \ref{BKH} 
the ${\cal Q}$-representation of $F$ injects into 
the ${\cal Q}$-representation of $E$.   
If the multiplicities of $S$ are all $1$ 
and $S'$  contains the source of  an arrow $\beta$  in $S$ then 
$S'$ contains $\beta$. }
\end{rem}

\begin{defin} We say that a  subgraph with multiplicities  of ${\cal Q} $ is  
a {\bf staircase} $S$  in a rectangle $R$ if all its multiplicities
 are $1$ and 
 the  graph of $S$ is a subgraph of  $R$  satisfying the following property:
if $V$ is a vertex of $S$ then the  arrows of $R$ having $V$ as 
source must be arrows of $S$ (and then also their sinks must be vertices
 of $S$).

We say that a  subgraph  with multiplicities  
of ${\cal Q}$  is a staircase if it is a staircase in 
some rectangle.
\end{defin}

Observe that a staircase $S$ in a rectangle $R$ has as matter of fact 
the form of a staircase with base and left side included respectively
in the base and the left side of $R$   as in the figure below. By Remark 
\ref{subquiv} the ${\cal Q}$-support of a    homogeneous subbundle
 of a homogeneous bundle whose
${\cal Q}$-support is a  rectangle is a staircase in the rectangle.

\begin{center}
\includegraphics[scale=0.26]{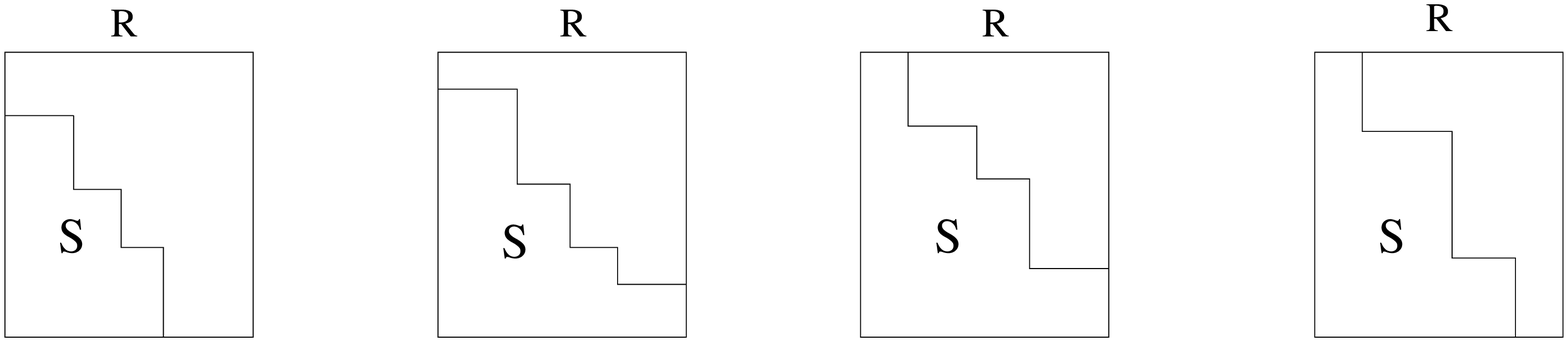}
\end{center}

\begin{notat} \label{scalini} (See the figure below). 
Given a staircase $S$ in a rectangle we define ${\cal V}_S$ to
 be the set of the vertices of $S$ that are not sinks of any arrow of $S$.
We call the elements of ${\cal V}_S$ the {\bf vertices of the steps}.
Let  ${\cal V}_S = \{V_1,..., V_k \} $  ordered in such a way the projection 
of $V_{i+1}$ on the base of  $R$ is on the left of the projection of
 $V_i$ $\forall i = 1,..., k-1$.

We define $R_i$ as the rectangles  with the right higher vertex equal 
to $V_i$ and left lower vertex equal to the left lower vertex of $R$.

For any $i = 1,...., k$,
we define the {\bf $i$-th horizontal step} $$H_i = R_i -R_{i-1}$$  ($R_0 = 
\emptyset$) and 
 the {\bf $i$-th vertical step} $$E_i = R_i -R_{i+1}$$  ($R_{k+1} = 
\emptyset$).
We define the {\bf $i$-th sticking out part} as $O_i = H_i \cap E_i $.

Let $s_i$ and $r_i$ be the lines containing respectively
the higher side of $O_i$ and  the right side of $O_i$.

For $i=1,..., k-1$,
let $S_i $ be the rectangle whose sides are on $r_i$, $s_{i+1}$, the line 
of the  base of the staircase and the line on the left side of the staircase.

Let $A_i $  $i=1,..., k-1$ and  $B_i $  $i=2,..., k$  be the rectangles 
$$ A_i = S_i -R_i -R_{i+1}
 \;\;\;\;\;\;\;\;\;\;\;\ B_i = S_{i-1} -R_i -R_{i-1}$$

\begin{center}
\includegraphics[scale=0.25]{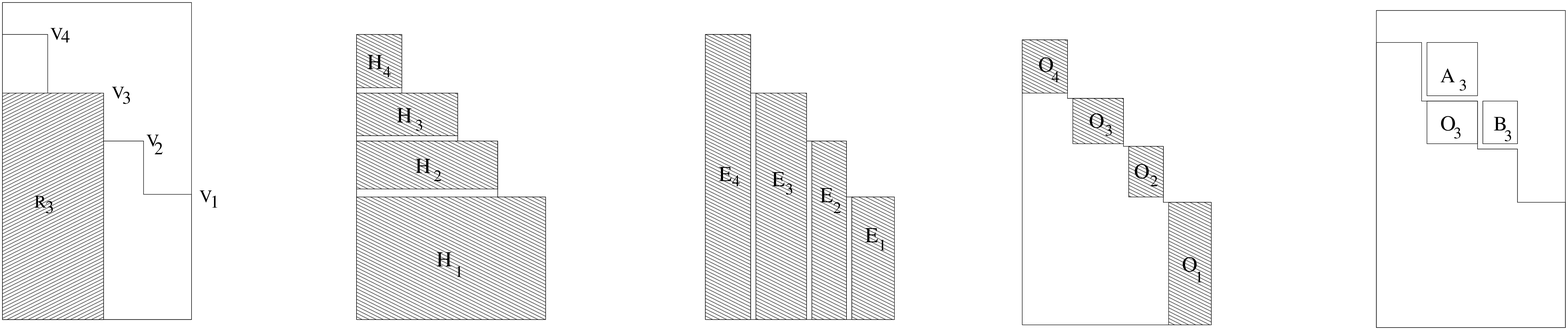}
\end{center}

\end{notat}

\section{Results on stability  and simplicity of  
 elementary  homogeneous bundles}

\begin{defin}
We say that a $G$-homogeneous  bundle is multistable if it is the tensor 
product of a stable $G$-homogeneous bundle 
and an irreducible $G$-representation. 
\end{defin}

\begin{theorem} \label{RF} {\bf (Rohmfeld, Faini)} 
i)  \cite{Rohm}  A homogeneous bundle
$E$ is semistable if and only if 
 $\mu (F) \leq \mu(E)$ for any subbundle $F$ of $E$  induced by a 
subrepresentation  of the $P$-representation  inducing  $E$.

ii) \cite{Fa}  A homogeneous bundle
$E$ is multistable if and only if 
 $\mu (F) < \mu(E)$ for any subbundle $F$ of $E$  induced by a 
subrepresentation  of the $P$-representation  inducing  $E$.

\end{theorem}

\begin{theorem} \label{stable} Let $ p,q \in {\bf N}$ with $p \geq q$ and 
$s >0$.
Let $E$ be the homogeneous vector bundle on ${\bf P}^2 = {\bf P}(V)$ defined
by the following exact sequence: \begin{equation} \label{as}
0 \rightarrow S^{p,q} V (-s)
\stackrel{\varphi}{\rightarrow}
 S^{p+s,q} V \rightarrow E \rightarrow 0 \end{equation} 
where $\varphi$ is a nonzero  $SL(V)$-invariant map.   Then $E$ is stable 
(in particular it is simple).
\end{theorem}

{\it Proof.} To show that $E$ is stable it is sufficient to show that it is 
multistable; in fact 
if $E$ is  the tensor product of a stable homogeneous 
vector bundle $E'$  with an $SL(V)$-representation $W$, then the minimal 
resolution of $E$ must be the tensor product of the minimal resolution of
 $E'$ with $W$ and from (\ref{as}) 
we must have 
$W ={\bf C}$.

To show that $E$ is multistable we consider the ${\cal Q}$-representation 
associated to $E$. In the figure we  show the sides of
the ${\cal Q}$-supports  
 of $S^{p,q}V(-s) $ and $S^{p+s,q}V$;
these supports  are rectangles  
(see Remark \ref{Spq}); thus the ${\cal Q}$-support 
of $E$  is the rectangle $R$:

\begin{center}
\includegraphics[scale=0.4]{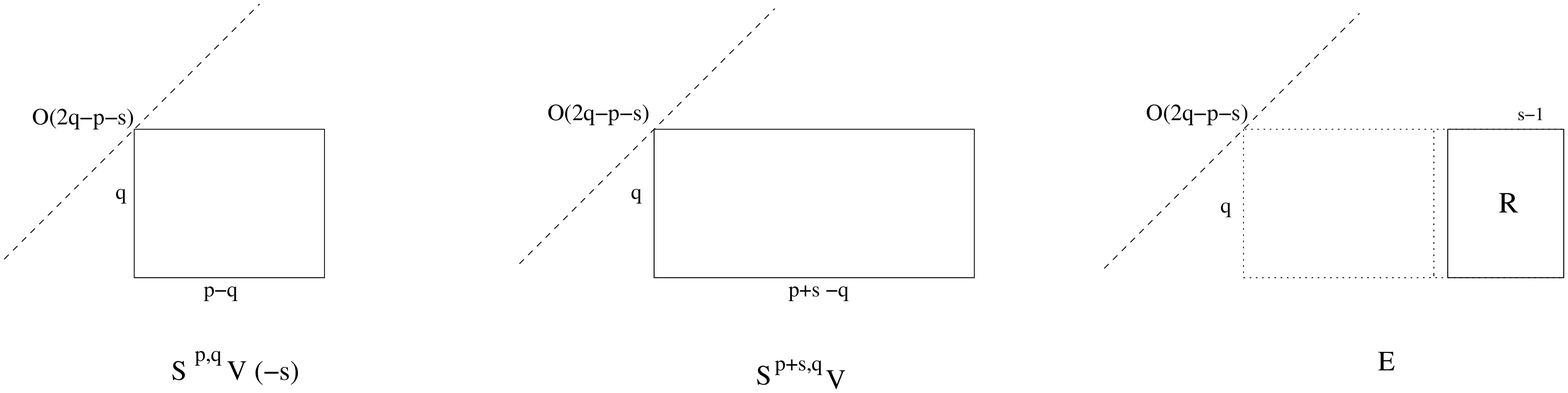}
\end{center}

By Theorem \ref{RF}, $E$ is multistable if
 $\mu (F) < \mu(E)$ for any subbundle $F$ of $E$  induced by a 
subrepresentation  of the $P$-representation  inducing  $E$.
  Observe that, by Remark \ref{subquiv},
 the support of the ${\cal Q}$-representation  
of any such subbundle $F$ must  be a staircase  $C$ in 
 $R$ 
and vice versa any ${\cal Q}$-representation  whose support is a 
 staircase $C$ in $R$
 is the ${\cal Q}$-representation of a
subbundle  $F$  of $E$ induced by a subrepresentation  of the 
$P$-representation inducing  $E$.

We will show by induction on the number $k$ of steps of $C$ 
that $\mu(C) < \mu(R)$ for any $C$ staircase in $R$.

\underline{$k=1$ }
In this case $C$  is a subrectangle in the  rectangle $R$.
 Thus this case follows from Lemma \ref{rettangoli}.

\smallskip

\underline{$k-1 \Rightarrow k$} 
We will show that, given a staircase $C$  in $R$  with $k$ steps, there 
exists a  staircase  $C'$ in $R$ with $k-1$ steps s.t. $\mu (C) \leq 
\mu(C') $. If we prove this, we conclude because $\mu(C)  \leq \mu (C') <
\mu(R) $, where the last inequality holds by  induction  hypothesis.

Let $C_1$ and $C_2$ be two staircases as in the figure:

\begin{center}
\includegraphics[scale=0.25]{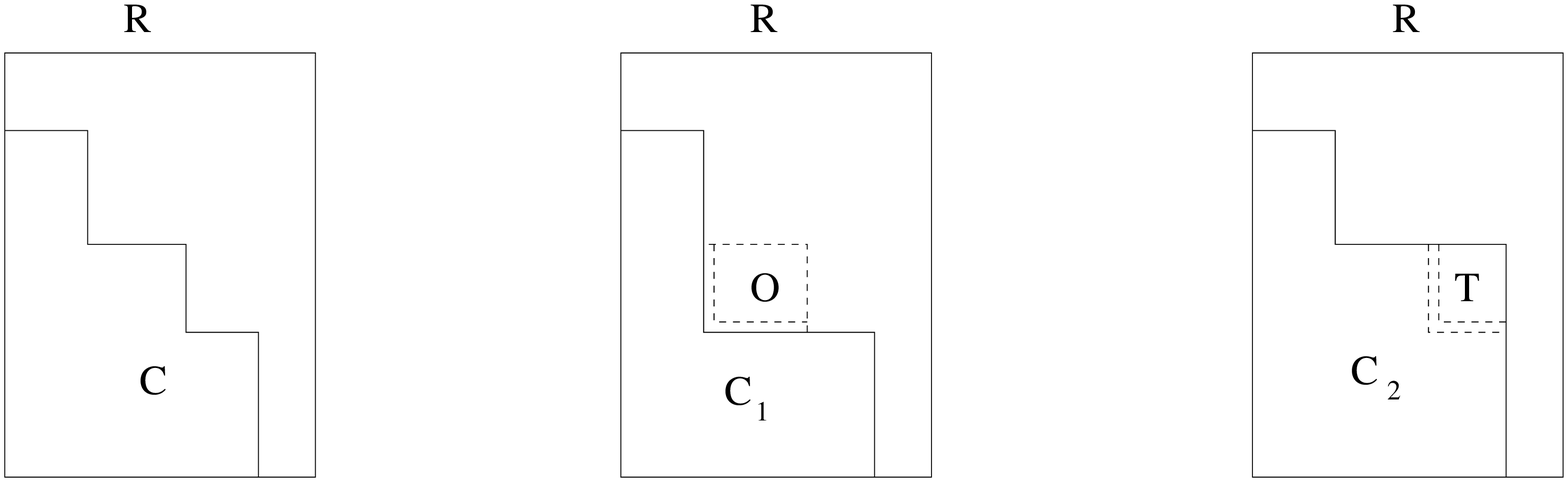}
\end{center}

that is $C_1$ and $C_2$  are staircases with $k-1$ steps 
obtained  from $C$ 
respectively ``removing and adding'' two rectangles $O$ and $T$.
Precisely $O$ is  a sticking out part $ O_i$ of $C$ for some $i$
 and    $T$ is a nonempty rectangle among the two rectangles
$A_i$, $B_i$ (see Notation \ref{scalini}).

If $\mu (C_1)  \geq \mu (C) $ we conclude at once.

Thus we can suppose that $\mu(C_1 ) < \mu(C)$. We state that in this case 
$\mu (C_2 ) \geq  \mu (C)$.   In fact: let $\mu (C_1)= \frac{a}{b} $, 
$\mu(O) =\frac{c}{d}$ and $ \mu(T) =\frac{e}{f}, $ where the numerators 
are the first Chern classes and the denominators the ranks;  
since $\mu (C_1)  < \mu (C) $, we have $ \frac{a}{b} < \frac{a+c}{b+d}$, 
thus $ \frac{a}{b}  < \frac{c}{d}$; besides by Lemma \ref{rettangoli}
 $\mu(O) < \mu(T) $, 
i.e. $\frac{c}{d} < \frac{e}{f}$; thus 
$\frac{a+c+e}{b+d+ f} \geq  \frac{a+c}{b+d}$ i.e. $\mu (C_2 ) \geq  \mu (C)$.  
\hfill  \framebox(7,7)

\medskip

Observe that, by Lemma \ref{soll} $i$, 
the simplicity statement of Theorem \ref{stable} is  equivalent to: 

\begin{cor} \label{lemmone}
 Let $ p,q ,s \in {\bf N}$ with $p \geq q$, $s > 0$ and 
let $A$ and $B$ be two linear maps s.t. the following diagram
commutes: $$\begin{array}{ccc} 
 S^{p,q } V \otimes S^s V & \stackrel{\pi}{\longrightarrow}
 & S^{p+s,q} V  \\ {\scriptsize A \otimes I}
\downarrow \;\;\;\;\;\;\; & & \;\;\; \downarrow {\scriptsize B}
\\S^{p,q }V \otimes S^s V& \stackrel{\pi}{\longrightarrow} & S^{p+s,q}V
\end{array}$$ where $\pi$ is  a nonzero
$SL(V)$-invariant projection (it is unique up
to multiple); then $A= \lambda I$ and $ B = \lambda I$ for some
$\lambda \in {\bf C}$.
\end{cor}

Now we want to prove Theorem \ref{pinco}; first it is necessary 
to prove several lemmas.

\begin{lemma} \label{lemmaL}
Let $p,q,s \in {\bf N}$  with $p\geq q$.

For every $ M \subset \{(s_1, s_2,s_3) \in {\bf N}^3
 |\; s_1 +s_2 +s_3 =s, \; s_2 \leq 
p-q,\;  s_3 \leq q \}$,
let ${\cal P}_M$ be the following statement: 
 for every  $V$  complex vector space of dimension $3$, the 
commutativity of the diagram  of bundles
on ${\bf P}(V)$: {\small $$\begin{array}{ccc}  S^{p,q } V (-s) &
\stackrel{\varphi}{\longrightarrow} & \oplus_{(s_1,s_2,s_3) \in M} S^{p+s_1,q
+s_2, s_3} V   \\ {\scriptsize A}
\downarrow \;\;\;\;\;\;\; & & \;\;\; \downarrow {\scriptsize B}
\\ S^{p,q }V (-s) & \stackrel{\varphi}{\longrightarrow} 
 & \oplus_{(s_1,s_2,s_3) \in M} 
S^{p+s_1,q +s_2, s_3} V  \end{array}$$ } \hspace*{-0.3cm} 
(where $A$ and $B$ are linear maps and the components of $\varphi$  are
nonzero $SL(V)$-invariant maps) implies  
$A= \lambda I$ and $B = \lambda I$ for some $\lambda \in {\bf C}$.

Let $$\{(s_1, s_2,s_3) \in {\bf N}^3 |\; s_1 +s_2 +s_3 =s, \; s_2 \leq p-q,\; 
s_3 \leq q \} =R \cup T$$ with $R \cap T = \emptyset $, $R \neq \emptyset$, 
$T \neq \emptyset $. 

Then ${\cal P}_R$ is true if and only if ${\cal P}_T$ is true.
\end{lemma}

{\it Proof.} Suppose ${\cal P}_{T}$ is true.  We want to show ${\cal P}_{R}$
 is true. 

Let $A$ and $B$ s.t. the diagram 
{\small $$\begin{array}{ccc}  S^{p,q } V (-s) &
\longrightarrow & \oplus_{(s_1,s_2,s_3) \in R} S^{p+s_1,q
+s_2, s_3} V    \\ {\scriptsize A}
\downarrow \;\;\;\;\;\;\; & & \;\;\; \downarrow {\scriptsize B}
 \\S^{p,q }V (-s) & \longrightarrow & \oplus_{(s_1,s_2,s_3) \in R} 
S^{p+s_1,q +s_2, s_3} V  \end{array}$$} \hspace{-0.26cm} 
 commutes. It is equivalent to the 
diagram 
{\small $$\begin{array}{ccc}  S^{p,q } V \otimes S^s V  &
\stackrel{\pi}{\longrightarrow} & \oplus_{(s_1,s_2,s_3) \in R} S^{p+s_1,q
+s_2, s_3} V    \\ {\scriptsize A \otimes I}
\downarrow \;\;\;\;\;\;\; & & \;\;\; \downarrow {\scriptsize B}
 \\S^{p,q }V \otimes S^s V 
& \stackrel{\pi}{\longrightarrow} & \oplus_{(s_1,s_2,s_3) \in R} 
S^{p+s_1,q +s_2, s_3} V  \end{array}$$} \hspace{-0.26cm} 
Thus $(A \otimes I) (Ker \pi) \subset 
Ker \pi$. Observe that $Ker \pi =  \oplus_{(s_1,s_2,s_3) \in T} S^{p+s_1,q
+s_2, s_3} V$. Then we get the following  commutative diagram: 
{\small  $$\begin{array}{ccc}  \oplus_{(s_1,s_2,s_3) \in T} S^{p+s_1,q
+s_2, s_3} V  & \longrightarrow & S^{p,q }V \otimes S^s V 
 \\ \downarrow \;\;\;\;\;\;\; & & \;\;\; \downarrow {\scriptsize A \otimes I}
 \\ \oplus_{(s_1,s_2,s_3) \in T} S^{p+s_1,q +s_2, s_3} V  & \longrightarrow & 
S^{p,q }V \otimes S^s V  \end{array}$$} \hspace{-0.26cm}
Let $W = V^{\vee}$. Substitute  $W^{\vee}$ for $V$
in the above diagram and dualize; the diagram  we obtain 
 is equivalent to the following commutative diagram of bundles on 
${\bf P}(W)$
{\small $$\begin{array}{ccc}  S^{p,q } W (-s) &
\longrightarrow & \oplus_{(s_1,s_2,s_3) \in T} S^{p+s_1,q
+s_2, s_3} W   \\ {\scriptsize A^{\vee} }
\downarrow \;\;\;\;\;\;\; & & \;\;\; \downarrow 
 \\S^{p,q }W (-s) & \longrightarrow & \oplus_{(s_1,s_2,s_3) \in T} 
S^{p+s_1,q +s_2, s_3} W  \end{array}$$} \hspace{-0.26cm}
By ${\cal P}_T$  we conclude that $A^{\vee} = \lambda I$ and then
 $A = \lambda I$.
\hfill  \framebox(7,7)

\begin{defin} We say that a  staircase is {\bf regular} if all the
 vertices of the steps (see Notation \ref{scalini})
 are on a line with angular  coefficient equal to $-1$.
\end{defin}

\begin{center}
\includegraphics[scale=0.28]{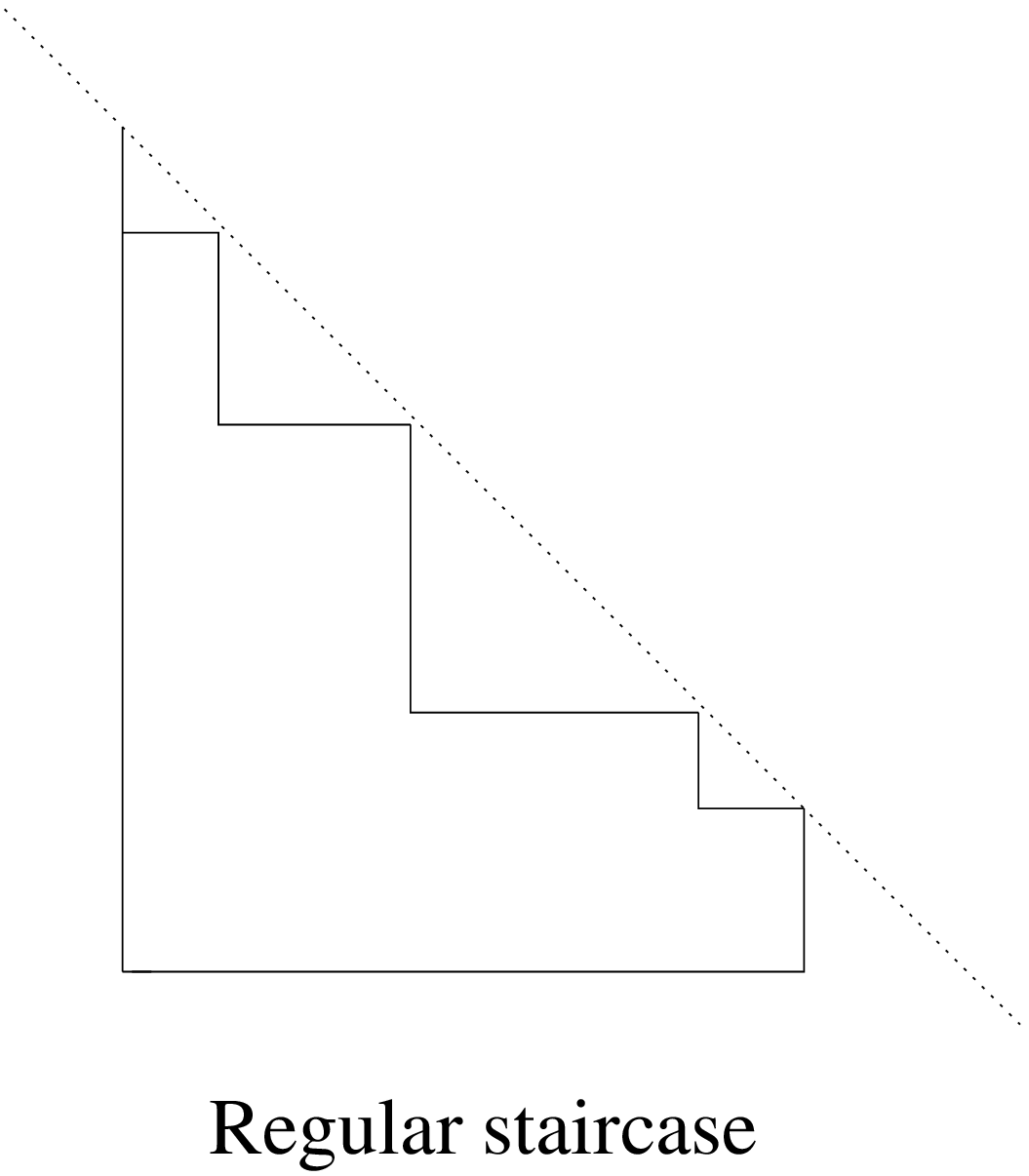}
\end{center}

\begin{lemma} \label{regstair}
The bundles whose support is a regular staircase are multistable.
\end{lemma}

{\it Proof.} 
{\sf Fact 1.} {\it  For any regular staircase we have 
$$\mu (H_i)  > \mu (H_{i-1}) \hspace{2cm}  \mu (E_i ) > \mu (E_{i+1})$$
for any $i$, where  $H_i$ are the   horizontal steps 
 and   $E_i$ are the vertical steps  (see Notation \ref{scalini}); recall 
that $H_{i+1}$ is  below $H_i$
and $E_{i+1}$ is  on left of $E_i$.}

{\sf Proof.} Obviously it is sufficient to prove the statement for a regular
staircase with two steps.
It is a freshman calculation  (even if a bit long)
 (use Lemma \ref{murettangolo}). 
 \hfill \framebox(4,4)

\smallskip

{\sf Fact 2.} {\it 
Let $S$ be a regular staircase. Then for every sticking out part $O$ of $S$
we have $$ \mu (O) > \mu (S-O)$$
Therefore } $$\mu (S) > \mu (S-O) $$
{\sf Proof.}
Let $b$ be the  line on which the base of $O$ is and  
let $l$ be the  line on which the left side of $O$ is. 
Let $T_1$ be the 
staircase  whose vertices are the vertices of $S$ that are  
either above $b$ or on $b$ and  on the left  of $l$ (see the figure below).
Let $T_2$ be the 
staircase  whose vertices are the vertices of $S$ that  are  
below  $b$ and either on the right of $l$  or on   $l$.

Let $K$ be the rectangle $$K= S -T_1 -T_2 -O$$

\begin{center}
\includegraphics[scale=0.3]{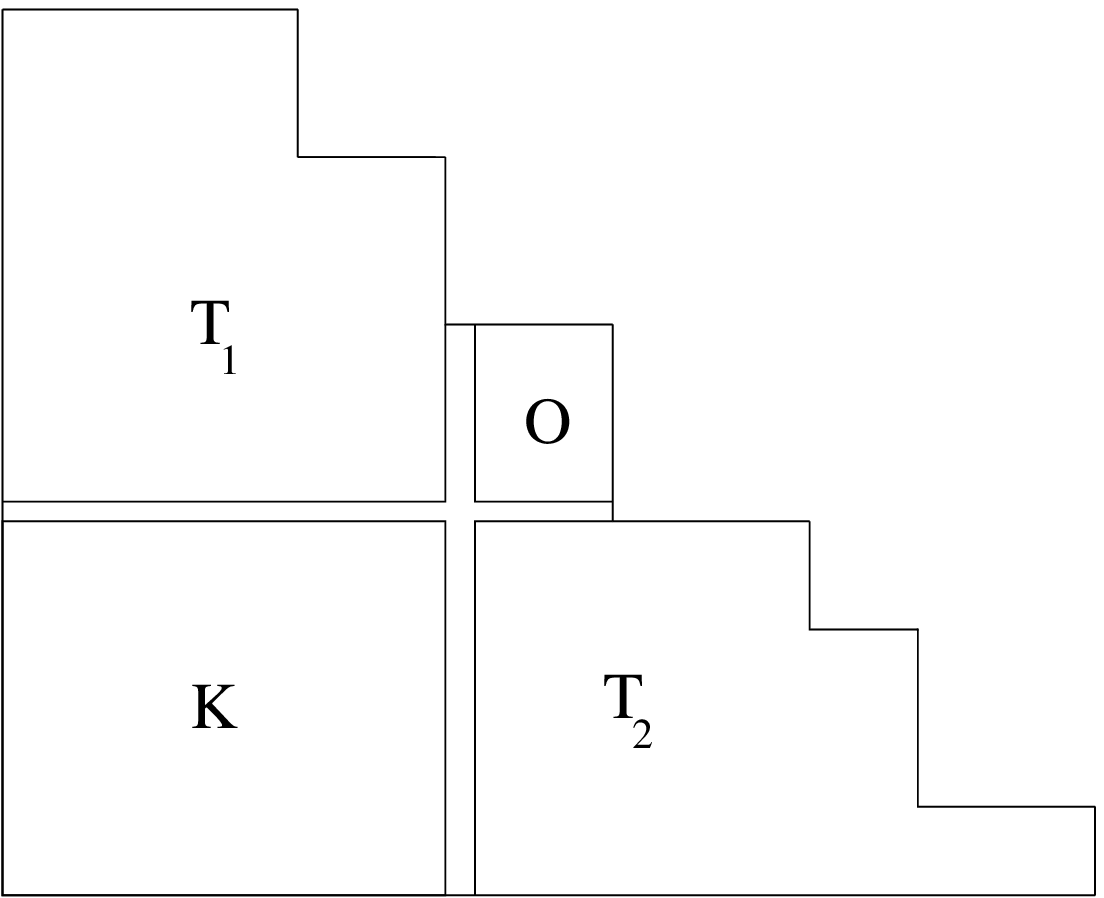}
\end{center}

By Lemma \ref{rettangoli} $\mu (O) > \mu (K) $.
Besides, by applying Fact 1 to the staircases $T_1 +O$ and $T_2 +O$
(where $T_i+O$ is the smallest staircase containing $T_i$ and $O$), we get 
 $$\mu (O) > \mu (T_1) \;\;\;\;\;\;\;\;\;\; \mu (O) > \mu (T_2)$$ 
Hence $\mu(O) > max \{\mu (K) , \mu(T_1) , \mu(T_2)\} \geq \mu(S-O)$ 
(see Remark \ref{gio}).
\hfill \framebox(4,4)

\smallskip

Now we are ready to prove that every bundle  s.t. its  ${\cal Q}$-support
  is a regular staircase $S$ is multistable. Let $C$ be the 
support 
of a ${\cal Q}$-representation subrepresentation of $S$ (thus again a 
staircase by Remark \ref{subquiv}). We want to prove $\mu(C) < \mu (S) $  by 
induction on the number $k$ of steps of $C$.

\underline{$k=1$}. The statement follows from
Lemma \ref{rettangoli} and Fact  1.

\underline{$k-1 \Rightarrow k$.} To prove this implication 
 we do induction on $$- length
(bd(C) \cap bd (S)) $$ 
where $bd $ denotes the border and the border of a staircase
 is the border of the part of the plane inside the staircase. 

Let $C$ be a staircase with $k$ steps  support 
of a subrepresentation of $S$. Let $O_i$ be the $i$-th 
sticking out part of $C$. 

$\bullet$
If $\mu (C -O_i) \geq \mu (C) $ for some $i$ we 
conclude at   once because $C-O_i $ has $k-1$ steps; thus by induction 
assumption $\mu(S) > \mu (C-O_i)$ and then $\mu (S) > \mu (C)$.  

$\bullet$ 
Thus we can suppose  $\mu (C) > \mu (C-O_i) $ $\forall i$
i.e.  $\mu (O_i) > \mu (C) $ $\forall i$.

Let $A'_i$ be the biggest rectangle in $A_i \cap S$ with the lower side
 equal to the lower side of $A_i$ and 
let $B'_i$ be the biggest rectangle in $B_i \cap S$ with the left side
 equal to the left side of $B_i$ (see Notation 
\ref{scalini} for the definition 
of $A_i$ and $B_i$).
\begin{center}
\includegraphics[scale=0.28]{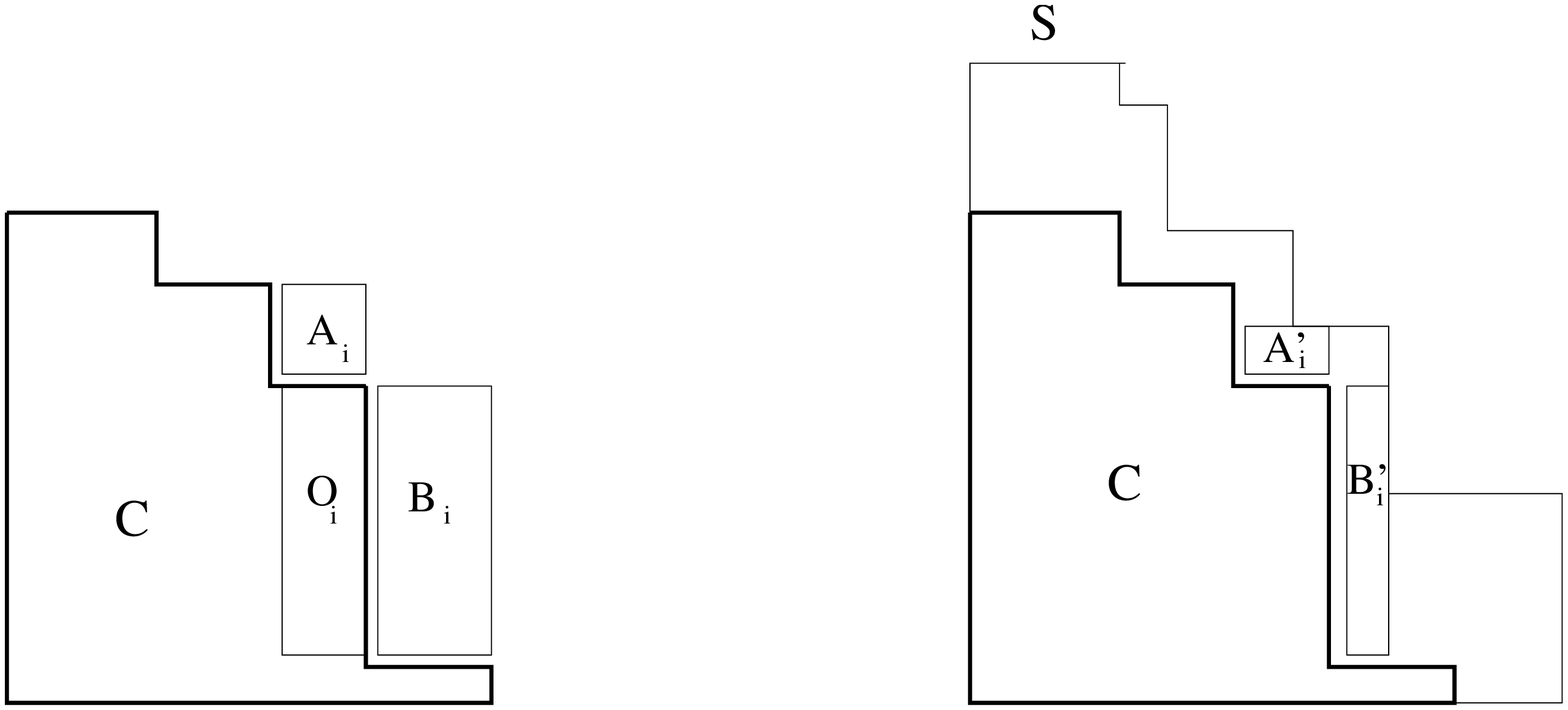}
\end{center}

$- $
Suppose $\exists i $ s.t. either $A'_i$ or $B'_i$ is not empty; for instance 
suppose $B'_i \neq \emptyset$. 

Since $\mu (B'_i) > \mu (O_i) $ by Lemma \ref{rettangoli}
 and $\mu(O_i)>  \mu (C) $ by assumption,
 we have $\mu (B'_i) > \mu (C)$ and thus
  \begin{equation} \label{st} \mu (C+ B'_i) > \mu (C) \end{equation} 
where $C+B'_i$ is the smallest staircase containing $C$ and $B'_i$.
If $B_i$ is a subgraph of $S$ i.e. $B_i= B'_i$, then $C + B'_i $ is a 
staircase  with $k-1$ steps thus, by induction assumption, $\mu (S) > \mu
 (C+ B'_i) $; hence 
 $\mu (S) > \mu (C)$  by (\ref{st}).

If $B_i$ is not a subgraph of $S$ i.e. $B_i \neq  B'_i$, then $length
 (bd (C+B'_i) \cap bd (S)  ) > length (bd(C) \cap bd (S))$ and by induction
 assumption $\mu (S) > \mu (C +B'_i) $; hence we conclude again
$\mu (S) > \mu (C)$  by (\ref{st}).

$-$  If $A'_i $ and $B'_i$ are empty  $\forall i$ then there
exists a chain of staircases 
 $C = S_0 \subset S_1 \subset ... \subset S_r =S$ s.t. $S_i $ is obtained 
from $S_{i+1}$ taking off one of its sticking out parts, thus, by Fact 2, we 
conclude.
\hfill  \framebox(7,7)

\begin{lemma} \label{Sltensor}
a) The ${\cal Q}$-support of $S^{l} Q (t) \otimes S^q V$ is (see also the
 figure below):

i) if  $l \geq q$,
 the subgraph of ${\cal Q}$, with all the  multiplicities  equal to $1$,
included  in an isosceles
right-angled  triangle with horizontal and vertical catheti of length $q$,
the direction of the hypotenuse  equal to NW-SE, the vertex opposite 
to the hypothenus equal to the lowest left vertex and 
equal to $S^l Q (t-q) $ 

ii)  if  $l < q$
 the subgraph of ${\cal Q}$, with all the multiplicities equal to $1$,
 included in a right-angled trapezium 
with horizontal bases, left side orthogonal to the bases, right side
with angular coefficient $-1$, length of the inferior base equal to $q$, 
the lowest left vertex equal to $S^{l} Q (t-q)$.

b) By duality (or directly) we can get an analogous statement for the
${\cal Q}$-support  of $S^{l} Q (t) \otimes S^{q, q}  V$ (see the figure).




\end{lemma}

{\it Proof.}
Use that $S^q V = \oplus_{i=0,...,q} S^{q-i} Q (-i)$ as $R$-representation  
and  Clebsch-Gordan's formula, see \cite{F-H}:  if $l \geq m $ 
$$S^l Q (t) \otimes S^m Q (r) = S^{l+m} Q  (t+r) \oplus S^{l+m-2} Q (t+r+1) 
\oplus  .... \oplus S^{l-m} Q (t +r +m)  $$ 
\hfill \framebox(7,7)

\begin{defin} We say that a staircase is {\bf completely 
regular} if it is equal to one of the  subgraphs of ${\cal Q}$ described in
a of Lemma \ref{Sltensor}.
\end{defin}

\begin{center}
\includegraphics[scale=0.28]{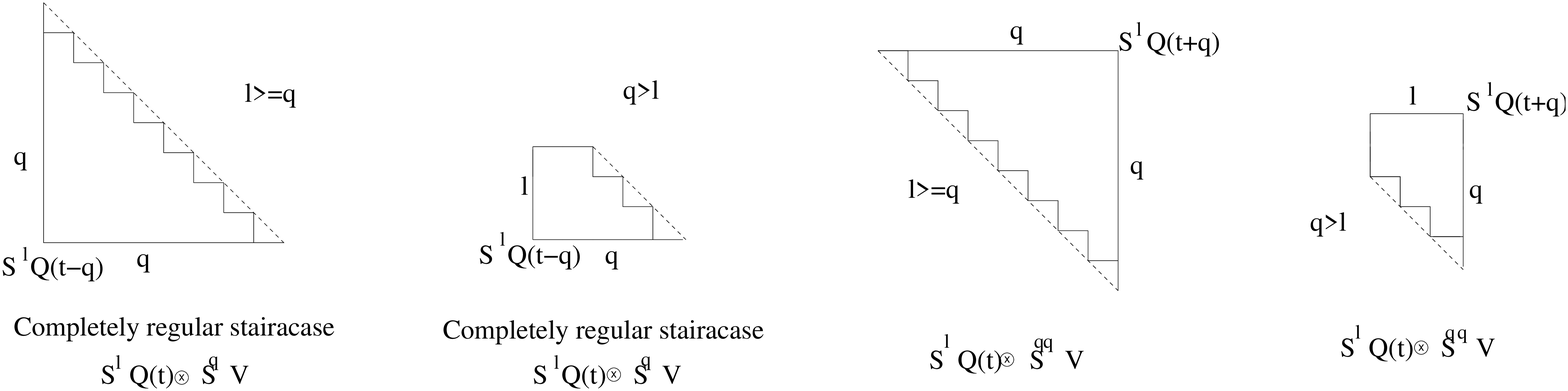}
\end{center}

\begin{lemma} \label{compregstair}
A bundle whose support is a regular staircase 
is stable if and only if the staircase is not completely regular.
\end{lemma}

{\it Proof.}
Observe that the ${\cal Q}$-support of 
$S^l Q (t) \otimes  S^{p,q} V$ has some  multiplicity
 $\geq 2$ if $p\neq q$ and 
$q\neq 0$; in fact
among  the vertices of the ${\cal Q}$-support of  $S^{p,q} V$ 
 there are $ {\cal O} (2q -p)$ and $ S^2 Q (-1 +2q -p)$, thus, by 
Clebsch-Gordan's formula,  $S^l Q 
(t +2q -p)$ occurs at least twice in $S^l Q (t) \otimes  S^{p,q} V$. 

By Lemma \ref{regstair} if a bundle $E$ has a regular staircase
as ${\cal Q}$-support, then $E = E' \otimes T$ where $E'$ is a 
stable   vector bundle and $ T$ is a vector space $SL(V)$-representation; 
if  $T= S^{p,q} V$ with $p \neq q $ and $ q \neq 0$ then by the previous remark
the ${\cal Q}$-support 
of $E' \otimes T$ has some  multiplicity $\geq 2$. Thus we must have
$T=S^{q,q} V$ or $T = S^{q} V$ and since the ${\cal Q}$-support 
of $E' \otimes T$  is a staircase only the last case is possible by Lemma 
\ref{Sltensor}.  
Thus we conclude by Lemma \ref{Sltensor}, in fact a regular staircase can be 
the disjoint union of $k$ completely regular staircases with the same
 length of  the base
if and only if  $k=1$ i.e. it is a completely regular stairacase (we 
can see this arguing on the upper part of the border).
\hfill  \framebox(7,7)

\begin{rem} \label{cofang}
{\rm If $p+q =p'+q'$, then the rectangles that are the 
${\cal Q}$-supports  of $S^{p,q} V (t) $ and  $S^{p',q'} V (t) $ 
have the highest right vertices on a line with angular coefficient $-1$.

If $p'' + q''  -t'' = p+q -t$ and $t'' <t$,
then the   highest right vertex of   the ${\cal Q}$-support 
of $S^{p'',q''} V (t'') $ is below this line.}
\end{rem}

\begin{lemma} \label{lemmaI}
Let  $p,q,s \in {\bf N}$ with $p\geq q$ and $s > 0$.

Let $T \subset \{(s_1, s_2,s_3) |\; s_1 +s_2 +s_3 =s, \; s_2 \leq p-q,\; 
s_3 \leq q \} $, $T \not \ni (s,0,0)$.
Let $A$ and $B$ be two linear maps s.t. the following diagram of bundles
on ${\bf P}(V) $ commutes: {\small $$\begin{array}{ccc}  S^{p,q } V (-s) &
\stackrel{\varphi}{\longrightarrow} & \oplus_{(s_1,s_2,s_3) \in T} S^{p+s_1,q
+s_2, s_3} V    \\ {\scriptsize A }
\downarrow \;\;\;\;\;\;\; & & \;\;\; \downarrow {\scriptsize B}
\\S^{p,q }V (-s) & \stackrel{\varphi}{\longrightarrow} 
& \oplus_{(s_1,s_2,s_3) \in T} 
S^{p+s_1,q +s_2, s_3} V  \end{array}$$} \hspace{-0.3cm}  where $\varphi $ is 
 an  $SL(V)$-invariant map with all its components nonzero.
 Then $A= \lambda I$ and $ B =
\lambda I$ for some $\lambda \in {\bf C}$.
\end{lemma}

{\it Proof.} Observe that 
  $A  (Ker (\varphi)) \subset  Ker (\varphi)$.
Thus we have a commutative diagram {\small
 \begin{equation}\label{ast1} \begin{array}{ccc}  Ker (\varphi) &
\longrightarrow &  S^{p,q} V (-s) \\ \;\;\;\;\;\;\;\;\;\;\;\downarrow
A|_{Ker(\varphi)}  & &  \downarrow A 
\\  Ker (\varphi) & \longrightarrow & S^{p,q} V (-s) \end{array}
\end{equation}} \hspace*{-0.3cm} 
Let $$ 0 \rightarrow R \rightarrow S \rightarrow Ker (\varphi) \rightarrow 
0 $$ be a minimal free resolution of $Ker(\varphi)$. 
By Lemma \ref{soll} 
the map $A|_{Ker(\varphi)} : Ker (\varphi) \rightarrow Ker (\varphi)$ induces 
a commutative diagram {\small 
$$\begin{array}{ccccccc} 0 \rightarrow  & R & \rightarrow &
 S & \rightarrow  & Ker(\varphi)  & \rightarrow 0
\\ & \downarrow  & & \downarrow  & & \;\;\;\;
\downarrow_{ A|_{Ker(\varphi)}}
 & \\  0 \rightarrow
& R & \rightarrow & S & \rightarrow  & Ker (\varphi)   & \rightarrow 0
\end{array}$$ } \hspace*{-0.3cm}
 Let $S_{max}$ be 
the direct sum of the summands of $S$ with maximum twist and let 
$S = S_{max} \oplus S'$; thus the previous diagram is      
{\small $$\begin{array}{ccccccccc} 0 \rightarrow  & R & \rightarrow &
 S_{max} & \oplus  & S' & \rightarrow  & Ker(\varphi)  & \rightarrow 0
\\ & \downarrow  & & \downarrow & \swarrow & \downarrow
 & & \;\;\;\; 
\downarrow_{ A|_{Ker(\varphi)}} & \\  0 \rightarrow
& R & \rightarrow & S_{max} &  \oplus &  S'& \rightarrow
  & Ker (\varphi)   & \rightarrow 0
\end{array}$$} \hspace{-0.3cm} Then we get a commutative diagram 
{\small $$ \begin{array}{ccc} S_{max}  & \rightarrow  & Ker(\varphi) 
\\  \downarrow  \alpha & & \;\;\;\;\downarrow A|_{Ker(\varphi)}  \\ 
S_{max}  & \rightarrow &  Ker (\varphi) 
\end{array}$$} \hspace*{-0.3cm}
 and, if $f$ is the composition of the map $S_{max } \rightarrow
 Ker (\varphi )$ with the  inclusion 
$ Ker (\varphi) \rightarrow S^{p,q} V (-s) $, by (\ref{ast1})
 we get the commutative diagram   
{\small \begin{equation} \label{ast2}
\begin{array}{ccc} S_{max}  &  \stackrel{f}{\rightarrow}  & S^{p,q} V (-s)
\\  \downarrow \alpha  & & \downarrow A  \\ 
S_{max}  & \stackrel{f}{\rightarrow} &  S^{p,q} V (-s) 
\end{array}\end{equation}} \hspace*{-0.3cm} and thus a commutative diagram 
{\small \begin{equation}  \label{ast3} 
\begin{array}{ccccccc} 0 \rightarrow &  Ker (f) & \rightarrow & 
S_{max}  &  \stackrel{f}{\rightarrow}  & Im (f) & \rightarrow 0 
\\  & \downarrow & & \downarrow \alpha  & & \downarrow \gamma 
 & \\  0 \rightarrow
 &  Ker (f) & \rightarrow & 
S_{max}  &  \stackrel{f}{\rightarrow}  & Im (f) & \rightarrow 0  
\end{array}
\end{equation}} \hspace*{-0.3cm} where $  \gamma = 
A|_{Im (f)}$.  

Now we will prove that 
 if $Im (f) $ is simple then $A$ 
is a multiple of the identity.

Let $0 \rightarrow  K \rightarrow M \rightarrow (Im f )^{\vee} \rightarrow 0 $
 be  a minimal free resolution of $ (Im f )^{\vee}$; 
for any $\beta : M \rightarrow M$ induced by $\gamma^{\vee}$
we have the following commutative  diagram: 
{\small  $$ \begin{array}{ccc} M  &  \rightarrow  & (Im \,f )^{\vee}
\\  \downarrow \beta  & & \downarrow  \gamma^{\vee} \\ 
M & \rightarrow  &  (Im \,f)^{\vee} 
\end{array}$$} \hspace*{-0.3cm}
and then by (\ref{ast3})
{\small 
 $$ \begin{array}{ccc} M  &  \stackrel{r}{\rightarrow}  & S_{max}^{\vee}
\\  \downarrow \beta  & & \downarrow  \alpha^{\vee} \\ 
M & \stackrel{r}{\rightarrow}  &  S_{max}^{\vee} 
\end{array}$$} \hspace*{-0.3cm}
 In particular, since $(Im f)^{\vee}$ is simple, $\gamma^{\vee}$ is
 a multiple of the identity, thus $\beta $ can  be taken equal to a multiple of
 the identity. 

Observe that all the 
components of $r$ are nonzero (because the map $M \rightarrow (Im \; f)^{\vee}$
is surjective  and all the components of $ f^{\vee}  : (Im \; f)^{\vee}
\rightarrow S_{max}^{\vee}$ are nonzero, since no component of $ S_{max}$
is sent to $0$ by  $f$); besides,   up to twisting, we can suppose $S_{max}$
is a trivial bundle and then $H^0(r^{\vee})^{\vee}$ is a projection.
Hence, since  $\beta $ is a multiple of the identity, $\alpha^{\vee}$  
(and then  $\alpha$) is a multiple of the identity.

Thus, by (\ref{ast2}) also $A$ is a multiple of the identity as we wanted
(twist (\ref{ast2}) by $s$ and consider $H^0 (\cdot^{\vee})^{\vee}$ of 
every map of  the obtained diagram; in this way we get a diagram  of 
$SL(V)$-representations whose horizontal maps are surjective, since they are 
$SL(V)$-invariant and nonzero and $S^{p,q} V $ is irreducible).

Observe that, by Remark \ref{cofang},  the ${\cal Q}$-support 
 of $Im f $ is a regular staircase in the
   ${\cal Q}$-support   of $S^{p,q} V (-s)$, since all the summands of
 $S_{max}$ have the same twist.

Thus by Lemma \ref{compregstair} we can conclude at once if the 
${\cal Q}$-support of   $Im f$ is not a completely regular staircase. 

Therefore we can suppose  the ${\cal Q}$-support 
of  $Im f$ is a completely regular staircase.

Observe that, by Remark \ref{cofang}, 
if the ${\cal Q}$-support $Im f $ is a completely regular staircase 
then $S= S_{max}$ and thus $Im f =  Ker  \varphi$.

Consider the following  exact sequence \begin{equation}\label{ast4}
0 \rightarrow Ker \varphi \rightarrow 
S^{p,q} V (-s)  \rightarrow Im \varphi \rightarrow 0 \end{equation}
Up to dualizing  we can suppose  that also  the ${\cal Q}$-support of 
 $ Im \varphi $ is of the kind {\it b} of Lemma \ref{Sltensor}.

  Thus the unique remaining cases   are the cases in 
which the sequence (\ref{ast4}) twisted by $s$ is one of the following:
$$ 
 Case \;A \;\;\;\;\;\;\;
0 \rightarrow S^q Q (-q-1) \otimes S^{p-q-1} V    \rightarrow S^{p,q} V  
  \rightarrow S^{p-q} Q \otimes S^{q,q} V   \rightarrow 0 $$ 
 $$ Case \;B \;\;\;\;\;\;\;
0 \rightarrow S^q Q (-q) \otimes S^{p-q} V    \rightarrow S^{p,q} V  
  \rightarrow S^{p-q} Q (1) \otimes S^{q-1,q-1} V  \rightarrow 0 $$

\begin{center}
\includegraphics[scale=0.32]{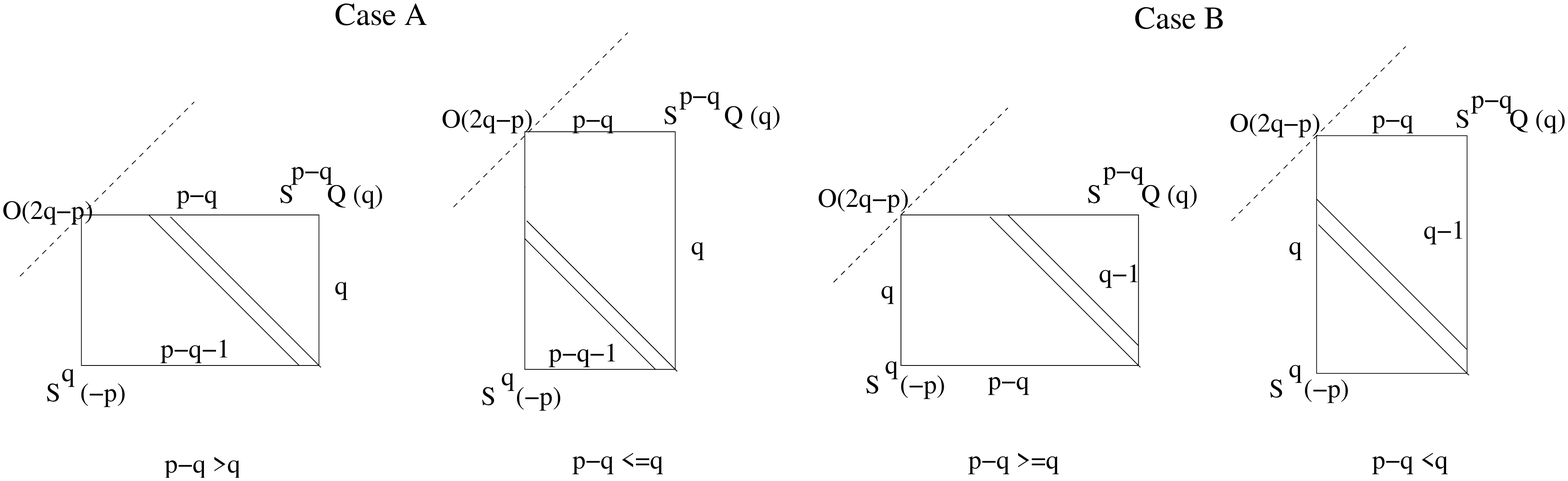}
\end{center}

Observe that the Case B is equivalent to Case A (by dualizing 
and considering $q'=p-q$). Thus it is sufficient to consider Case A.
By (\ref{ast1})  we get a commutative diagram 
{\small $$ \begin{array}{ccccccc}
0  \rightarrow 
 & S^q Q (-q-1) \otimes S^{p-q-1} V  &  \rightarrow & S^{p,q} V  
&  \rightarrow & S^{p-q} Q \otimes S^{q,q} V &  \rightarrow 0  \\
 & \downarrow & & \downarrow & & \downarrow & \\ 
  0 
\rightarrow & S^q Q (-q-1) \otimes S^{p-q-1} V  &  \rightarrow & S^{p,q} V  
 & \rightarrow & S^{p-q} Q \otimes S^{q,q} V  & \rightarrow 0 
\end{array} $$}
\hspace{-0.26cm}
By taking the cohomology (in particular $H^0 $) we get   
{\small $$ \begin{array}{ccccccc} 0 & \rightarrow &
0  & \rightarrow  &   S^{p,q} V  
 & \rightarrow & S^{q,q} V \otimes S^{p- q} V \\
   \downarrow  & & \downarrow & & \downarrow  \\ 
0 & \rightarrow  & 0 &  \rightarrow  & S^{p,q} V  
 & \rightarrow & S^{q,q} V \otimes S^{p-q} V \end{array} $$}
\hspace{-0.26cm}
since $H^0 (S^q Q (-q-1) ) =0$ and $ H^0 (S^{p-q} Q) = S^{p-q} V$
(to calculate  the cohomology of $S^l Q  (t) $ use for instance its
  minimal resolution  $ 0 \rightarrow S^{l-1} V (t-1) \rightarrow 
S^{l} V (t) \rightarrow S^l Q (t) \rightarrow 0$). We conclude by  applying 
Lemma \ref{lemmone}  to the dual of the right part of the diagram.
\hfill  \framebox(7,7)

\bigskip

From Lemmas \ref{lemmaI} and \ref{lemmaL} we deduce at once: 

\begin{cor} \label{lemmonetris} Let $ p,q ,s\in {\bf N}$ with $p\geq q$. 
Let $A$ and $B$ be two linear maps s.t. the following diagram of bundles
on ${\bf P}(V) $ 
commutes: {\small $$\begin{array}{ccc}  S^{p,q } V (-s) &
\stackrel{\varphi}{\longrightarrow} & W  \\ {\scriptsize A \otimes I}
\downarrow \;\;\;\;\;\;\; & & \;\;\; \downarrow {\scriptsize B}
\\ S^{p,q }V (-s) & \stackrel{\varphi}{\longrightarrow} & W
\end{array}$$} \hspace*{-0.3cm}
 where $W$ is  a non trivial  $SL(V)$-submodule of   
$ S^{p, q} V \otimes S^s V$ and
all the components of $\varphi$ are  nonzero $SL(V)$-invariant maps.
 Then $A= \lambda I$ and $ B =
\lambda I$ for some $\lambda \in {\bf C}$.
\end{cor}

\medskip

{\it Proof of Theorem \ref{pinco}.} The case $p=0$ is trivial. Thus we can
 suppose $p > 0$. 
 By Corollary \ref{lemmonetris},
the only thing we have to prove is that if in $W \otimes {\cal O}$
 there are two copies of an irreducible bundle $F$, i.e. 
$W = F \oplus F \oplus W'$, then $E$ is not simple: in fact 
the following diagram induces an automorphism on $E$ not multiple 
of the identity:
{\small  $$
\begin{array}{ccccccccccc} 0 \rightarrow  & S^{p,q} V (-s) & \rightarrow &
 F & \oplus & F & \oplus & W'  & \rightarrow  & E  & \rightarrow 0
\\ & \downarrow {\small I} & & -I \downarrow \;\;\;\;\;  
& \swarrow_{ 2{\small I}} &
\;\;\; \downarrow {\small I} &  & \;\; \downarrow {\small I}
 & & \downarrow & \\  0 \rightarrow
& S^{p,q} V (-s) & \rightarrow &  F & \oplus & F & \oplus & W' &
\rightarrow  & E  & \rightarrow 0
\end{array}$$}
\hspace{-0.26cm}
 (see Lemma \ref{soll}).
\hfill  \framebox(7,7)

\medskip

The following theorem gives a precise criterion to see when a regular 
elementary homogeneous bundle $E$ on ${\bf P}^2$
is stable or simple in the case the difference 
of the twists of the first bundle and of the middle bundle of the minimal 
free resolution of $E$ is $1$.

\begin{theorem} \label{1regular}
Let $E$ be a homogeneous bundle on ${\bf P}^2$
with minimal free resolution {\small
$$ 0 \rightarrow S^{p,q}V \stackrel{\varphi}{\rightarrow} 
\oplus_{\alpha\in {\cal A}}
S^{p+s_1^{\alpha},q+s_2^{\alpha},s_3^{\alpha}}V(1)
\rightarrow E \rightarrow 0 $$ } \hspace*{-0.3cm}
with $p \geq q \geq 0$, $s_i^{\alpha} \in {\bf N}$,
 $1=s_1^{\alpha}+s_2^{\alpha}+s_3^{\alpha}$, ${\cal A}$ finite subset of 
indices and all the components of $\varphi$ 
nonzero $SL(V)$-invariant maps. Then  

(i) $E$ is simple if and only if its minimal free resolution has one of 
the following forms:
$$0 \rightarrow S^{p,q}V \rightarrow
S^{p+1,q}V(1) \rightarrow E \rightarrow 0$$
$$0 \rightarrow S^{p,q}V \rightarrow
S^{p+1,q}V(1) \oplus S^{p,q+1}V(1) \rightarrow E \rightarrow 0$$
 with $q\neq 0$,
$$0 \rightarrow S^{p,q}V \rightarrow
S^{p+1,q}V(1) \oplus S^{p,q,1}V(1) \rightarrow E \rightarrow 0$$
 with $q \neq 0$ and $p \neq q$.

(ii) $E$ is stable if and only if $E$ is simple and, moreover, in the 
case its minimal resolution   is of the third type
 $2q \geq p >q$.
\end{theorem}

{\it Proof} 
(i) follows from Theorem \ref{pinco}. So it is enough to check
when the bundles described in (i) are stable.
The first case follows from Lemma \ref{lemmone}.
For the second case note that the ${\cal Q}$-support
  of $E$ is a not completely regular staircase, thus we conclude by Lemma
\ref{compregstair}.
In the third case the ${\cal Q}$-support  of $E$ is the following:
\begin{center}
\includegraphics[scale=0.5]{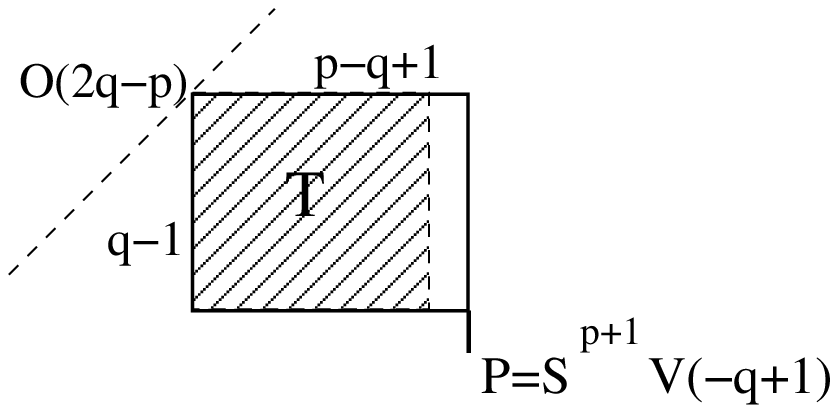}
\end{center}
The bundle $E$ is multistable if and only if $ \mu(T) <
 \mu(E) $ and $\mu (P) < \mu (E)$
and, by using  Lemma \ref{murettangolo}, one can show that this is true 
if and only if $ -5p +2 q- 2 -4 p^2 +3pq +p^2q -p^3 <0$ and $ 3 p^2 + p^3 -8
 pq 
-2 p^2 q -4 -8 q < 0$   
respectively. The first inequality always holds since $p\geq q$.
The last inequality is equivalent to $p\leq 2q$ 
since  its first member is equal to
$(p+2)^2 (p-1-2q) $. Thus $E$ is multistable if and only if $p \leq
2q$.   In this case, by Lemma \ref{Sltensor}, it is  stable if and only if 
 $ p \neq q$.   
\hfill \framebox(7,7)

\medskip

 In \cite{Fa} an example of a simple unstable homogeneous 
bundle on ${\bf P}^2$ is exhibited. Theorem \ref{1regular}
shows infinite examples of such  bundles.

We end with a theorem which studies  the simplicity of elementary 
homogeneous bundles.  

\begin{theorem} \label{fibomsem} Let $p \geq q $.
Let $E$ be the homogeneous vector bundle on ${\bf P}^2
= {\bf P}(V)$  defined
by the following exact sequence
(${\cal A}$ a finite
set of indices): $$0 \rightarrow S^{p,q} V
\stackrel{\varphi}{\rightarrow} 
\oplus_{\alpha \in {\cal A}} S^{p+s_1^{\alpha},q
+s_2^{\alpha}, s_3^{\alpha}} V (s^{\alpha}) \rightarrow E
\rightarrow 0$$ 
where the components $\varphi_{\alpha}$ of  $\varphi$ 
   are $SL(V)$-invariant maps and
 $s_1^{\alpha} + s_2^{\alpha} + s_3^{\alpha}
=s^{\alpha} $ $\forall \alpha \in {\cal A}$.  
Then $E$ is simple if and only if the following five conditions hold:

a) $\not \exists \alpha, \beta \in {\cal A}$ s.t. $s_i^{\alpha}
\leq s_i^{\beta}$ $i=1,2,3$

b) all  the components $\varphi_{\alpha} $  of   $\varphi$
are nonzero, in particular
 $p \geq q+ s_2^{\alpha} $  and $q \geq s_3^{\alpha}$
 $\forall \alpha \in {\cal A}$

c) $\forall \alpha , \beta \in {\cal A}$ s.t. $s^{\alpha} >
s^{\beta }$, we have that  if $s_2^{\beta } > 0 $ then $s_3^{\beta
} + s^{\alpha} - s^{\beta} < q+1  $ and if $s_1^{\beta } > 0 $
then $q+ s_2^{\beta } + s^{\alpha} - s^{\beta} < p + 1  $

d) $\forall \alpha , \beta, \gamma \in {\cal A}$ s.t. $s^{\alpha}
= s^{\beta }    < s^{\gamma} $, we have $ s^{\gamma} -
s^{\alpha} \geq  max \{|s_i^{\alpha }- s_i^{\beta }| \; i=1,2,3\}$

e) if $p > 0 $ and  $\exists \; s$ s.t. $s^{\alpha} = s$ 
 $\forall \alpha \in {\cal A}$ then 
$\oplus_{\alpha \in {\cal A}} S^{p+s_1^{\alpha},q
+s_2^{\alpha}, s_3^{\alpha}} V  \neq S^{p,q}V \otimes S^s V$.
\end{theorem}


{\it Sketch of the proof.}
Any automorphism  $\eta$ of $E$ induces maps $A$ and $B$ as in Lemma 
\ref{soll} and given $A$ and $B$ we have an automorphism of $E$.

We call $B_{{\cal K},{\cal J}} :\oplus_{\alpha \in {\cal
J}} S^{p+s_1^{\alpha},q +s_2^{\alpha}, s_3^{\alpha}} V
(s^{\alpha}) \rightarrow \oplus_{\alpha \in {\cal K}}
S^{p+s_1^{\alpha},q +s_2^{\alpha}, s_3^{\alpha}} V (s^{\alpha}) $
the map induced by $B$,  $\forall {\cal J}, {\cal K} $ of ${\cal A}$. 
We denote $B_J =B_{J,J}$, $B_{\alpha, \beta}=  B_{\{\alpha\},
 \{\beta\}} $ and $B_{\alpha}=  B_{\{\alpha\}} $ for short.
Besides  $ \varphi_{\cal J} : S^{p,q} V \rightarrow
\oplus_{\alpha \in {\cal J}} S^{p+s_1^{\alpha},q +s_2^{\alpha},
s_3^{\alpha}} V (s^{\alpha})$   denotes the map induced by $\varphi$,
 $\forall {\cal J} \subset {\cal A}$.

We show now that the simplicity of $E$ implies {\it a,b,c,d,e};
{\it a},{\it b} and {\it e} are  left to the reader. 



{\it c)} 
Observe that if for some $\overline{\alpha},
\overline{\beta} \in {\cal A}$ there exists a nonzero map $$\gamma:
S^{p+s_1^{\overline{\beta}},q +s_2^{\overline{\beta}},
s_3^{\overline{\beta}}} V (s^{\overline{\beta}}) \rightarrow
S^{p+s_1^{\overline{\alpha}},q +s_2^{\overline{\alpha}},
s_3^{\overline{\alpha}}} V (s^{\overline{\alpha}})$$ s.t. $\gamma
\circ \varphi_{\overline{\beta}} =0$ then $E$ is not simple
(take $B_{\overline{\alpha},\overline{\beta}}=
\gamma$, $A=I$, $B_{\alpha}=I$ $\forall \alpha \in {\cal A}$,
$B_{\alpha, \beta }=0$ $\forall (\alpha, \beta) \neq
(\overline{\alpha}, \overline{\beta})$ and use  Lemma \ref{soll}).
Such a $\gamma $ exists if and only if 
there exists a nonzero map $$   \Gamma:
S^{p+s_1^{\beta},q +s_2^{\beta}, s_3^{\beta}} V \otimes 
S^{s^{\overline{\alpha}} -s^{\beta}} V
\rightarrow
 S^{p+s_1^{\overline{\alpha}},q 
+s_2^{\overline{\alpha}}, s_3^{\overline{\alpha}}} V
$$ s.t. $ \Gamma \circ 
H^0((\varphi_{\beta}(-s^{\overline{\alpha}}))^{\vee})^{\vee}=0$
This is equivalent to the non surjectivity
of $H^0((\varphi_{\beta}(-s^{\overline{\alpha}})^{\vee})^{\vee}$
   (which is the injection followed by the 
projection $S^{p,q} V \otimes  S^{s^{\overline{\alpha}}} V \rightarrow
S^{p,q} V \otimes  S^{s^{\beta}} V \otimes   
S^{s^{\overline{\alpha}} -s^{\beta}} V  \rightarrow 
S^{p+s_1^{\beta},q +s_2^{\beta}, s_3^{\beta}} V \otimes
S^{s^{\overline{\alpha}} -s^{\beta}} V$)
and we conclude  by Pieri's formula.

{\it d)} 
Obviously if for some $ \alpha , \beta, \gamma \in {\cal A}$ s.t.
$s^{\alpha} = s^{\beta }    < s^{\gamma} $ there exist  nonzero
maps
$$ \begin{array}{c} \delta: S^{p+s_1^{\alpha},q
+s_2^{\alpha}, s_3^{\alpha}} V (s^{\alpha}) \rightarrow
S^{p+s_1^{\gamma},q +s_2^{\gamma}, s_3^{\gamma}} V (s^{\gamma}) \\
   \delta': S^{p+s_1^{\beta},q +s_2^{\beta}, s_3^{\beta}} V
(s^{\beta}) \rightarrow S^{p+s_1^{\gamma},q +s_2^{\gamma},
s_3^{\gamma}} V (s^{\gamma}) \end{array} $$
s.t. $\delta \circ \varphi_{\alpha} + \delta' \circ
\varphi_{\beta} =0$ then $E$ is not simple (take
$A=I$, $B_{\alpha} =I$  $\forall \alpha \in {\cal A}$,
$B_{\gamma,\alpha} =\delta $, $B_{\beta,\alpha} =\delta' $ and
$B_{\epsilon,\lambda} =0 $ $\forall (\epsilon, \lambda) \neq
(\gamma,\alpha), (\gamma,\beta) $).
Such $\delta$ and $ \delta'$ exist
 if and only if 
{\small $$
 S^{p,q} V \otimes S^{s^{\gamma}} V   \stackrel{
 H^0(\varphi_{\alpha}(-s^\gamma)^{\vee})^{\vee} \times
 H^0(\varphi_{\beta}(-s^\gamma)^{\vee})^{\vee}}{ \longrightarrow} 
(S^{p+s_1^{\alpha},q +s_2^{\alpha}, s_3^{\alpha}} V \otimes
S^{s^{\gamma} -s^{\alpha}} V) \oplus (S^{p+s_1^{\beta},q +s_2^{\beta},
s_3^{\beta}} V \otimes S^{s^{\gamma} -s^{\alpha}}V)$$} 
\hspace{-0.1cm}is not surjective; 
since  $H^0(\varphi_{\alpha}(-s^\gamma)^{\vee})^{\vee} $ and $
 H^0(\varphi_{\beta}(-s^\gamma)^{\vee})^{\vee}$
 are surjective by {\it c},
this is true if and only if 
 $S^{p+s_1^{\alpha},q +s_2^{\alpha}, s_3^{\alpha}} V \otimes
S^{s^{\gamma} -s^{\alpha}} V $ and $S^{p+s_1^{\beta},q +s_2^{\beta},
s_3^{\beta}} V \otimes S^{s^{\gamma} -s^{\alpha}} V$ have a nonzero
subrepresentation in common 
and we conclude.


Suppose now {\it a,b,c,d,e} hold.
We may assume $p >0$.
Let $\eta $ be an automorphism of $E$ and let $B$ and $A$ be the induced maps
as above.

Let $ {\cal A}  =  {\cal A}' \cup  {\cal A}'' \cup ...$ be disjoint
union s.t. $  {\cal A}' $ is the set of indices $\alpha $ in $ {\cal
A} $ s.t. $s_{\alpha}$ is the minimum of $\{s_{\alpha}|\; \alpha
\in {\cal A} \}$, $  {\cal A}'' $ is the set of indices $\alpha $ in
$ {\cal A} - {\cal A}'$ s.t. $s_{\alpha}$ is the minimum of
$\{s_{\alpha}|\; \alpha \in {\cal A}-{\cal A}' \}$ and so on.
Let $s' = s^{\alpha}$ for $\alpha \in {\cal A}'$ and $s'' =
s^{\alpha}$ for $\alpha \in {\cal A}''$ and so on. 

We have $\varphi_{{\cal A}'} \circ A = B_{{\cal A}'} \circ \varphi_{{\cal
A}'}$; thus, by Corollary 
\ref{lemmonetris},    $A= \lambda I$ and $B_{{\cal A}'} =\lambda I$.

Besides  we have $$
\varphi_{{\cal A}''} \circ A = B_{{\cal A}''} \circ \varphi_{{\cal
A}''} +  B_{ {\cal A}'', {\cal A}'} \circ \varphi_{{\cal A}'} $$
thus we get $ (\lambda I-B_{{\cal A}''}) \circ \varphi_{{\cal
A}''} =  B_{ {\cal A}'', {\cal A}'} \circ \varphi_{{\cal A}'} $.
By applying $H^{0}(\; \cdot^{\vee})^{\vee}$,  we obtain
 {\small   $$ \begin{array}{cccc}S^{p,q} V
\otimes S^{s''} V  & \stackrel{H^0 (\varphi_{{\cal A}''}^{\vee})^{\vee}}
 {\rightarrow} 
& \oplus_{\alpha \in {\cal A}''}
S^{p+s_1^{\alpha},q +s_2^{\alpha}, s_3^{\alpha}} V
\\ H^0 (\varphi_{{\cal A}'}^{\vee})^{\vee}
 \downarrow &
  &\downarrow  H^0 ((\lambda I -B_{{\cal A}''})^{\vee})^{\vee}
 \\  \oplus_{\alpha \in {\cal A}'}
S^{p+s_1^{\alpha},q +s_2^{\alpha}, s_3^{\alpha}} V \otimes S^{s''
-s' } V & \stackrel{H^0 (B_{{\cal A}'',{\cal  A}'}^{\vee})^{\vee}
 }{\rightarrow} 
 & \oplus_{\alpha \in {\cal A}''}
S^{p+s_1^{\alpha},q +s_2^{\alpha}, s_3^{\alpha}} V \end{array} $$}
\hspace*{-0.3cm}
  By {\it a} we have that
$\oplus_{\alpha \in {\cal A}''} S^{p+s_1^{\alpha},q +s_2^{\alpha},
s_3^{\alpha}} V$ is in the kernel of the map $$H^0 (\varphi_{A'}^{\vee} 
)^{\vee} : S^{p,q} V \otimes
S^{s''} V \rightarrow \oplus_{\alpha \in {\cal A}'} S^{p+
s_1^{\alpha},q +s_2^{\alpha}, s_3^{\alpha}} V \otimes S^{s''
-s'} V$$ thus $ \lambda I - B_{{\cal A}''} =0$, then
 $ B_{{\cal A}''} = \lambda I$ and  $B_{ {\cal A}'', {\cal
A}'} \circ \varphi_{{\cal A}'} =0  $. Hence  by  {\it c} and {\it d}
we get $B_{ {\cal A}'', {\cal A}'} =0 $ (arguing as in 
the proof of the other implication).

By induction on the number of the subsets ${\cal A}', {\cal
A}'',....$  we conclude.
  \hfill  \framebox(7,7)

\end{document}